\documentclass{svjour3}       
\RequirePackage{fix-cm}
\makeatletter
\let\cl@chapter\undefined
\makeatletter
\usepackage{amsmath,amsfonts}
\usepackage[linesnumbered,ruled,commentsnumbered]{algorithm2e}
\usepackage{algorithmicx}
\usepackage{array}
\usepackage[caption=false,font=normalsize,labelfont=sf,textfont=sf]{subfig}
\usepackage{textcomp}
\usepackage{stfloats}
\usepackage{booktabs}
\usepackage{multirow}
\usepackage{hyperref}
\usepackage{nicematrix}
\usepackage{cleveref}
\usepackage{url}
\usepackage{tikz}
\usepackage{verbatim}
\usepackage{graphicx}
\usepackage{pgfplots}
\usetikzlibrary{shapes.geometric}
\usepgfplotslibrary{groupplots}
\usepackage{pgfplotstable}
\pgfplotsset{compat=newest}
\usetikzlibrary{arrows}
\usetikzlibrary{decorations}
\usetikzlibrary{decorations.markings}
\usetikzlibrary{shapes}
\usetikzlibrary{intersections}
\usetikzlibrary{patterns}
\usetikzlibrary{pgfplots.statistics}
\Crefname{figure}{Fig.}{Figs.}
\usepackage{adjustbox}
\usetikzlibrary{shapes.geometric, arrows, fit}
\usetikzlibrary{automata,positioning}
\usetikzlibrary{shapes.multipart}
\tikzset{
	startstop/.style =
	{rectangle, rounded corners, minimum width=2cm, minimum height=0.5cm,text centered, draw=black, fill=red!30, line width=1pt},
	io/.style = 
	{trapezium, trapezium left angle=70, trapezium right angle=110, minimum width=3cm, minimum height=1cm, text centered, text width=3cm,draw=black, fill=blue!30},
	process/.style = {rectangle, rounded corners, minimum width=3cm, minimum height=1cm, text centered, text width=3cm, draw=black, fill=orange!30},
	dimensss/.style = {rectangle, minimum width=1cm, minimum height=0.75cm, text centered, text width=3.25cm, draw=black, fill=red!20},
	scheme/.style = {rectangle, minimum width=1cm, minimum height=0.5cm, text centered, text width=1cm, draw=black, fill=blue!30},
	dimens/.style = {rectangle, minimum width=1cm, minimum height=0.75cm, text centered, text width=3.25cm, draw=black, fill=blue!30},
	selegh/.style = {rectangle, minimum width=1cm, minimum height=0.75cm, text centered, text width=5cm, draw=black, fill=blue!30},
	select/.style = {rectangle, minimum width=2cm, minimum height=0.75cm, text centered, text width=2.25cm, draw=black, fill=blue!30},
	diamet/.style = {rectangle, minimum width=1cm, minimum height=0.75cm, text centered, text width=2cm, draw=black, fill=blue!30},
	diamest/.style = {rectangle, minimum width=1cm, minimum height=0.75cm, text centered, text width=1.75cm, draw=black, fill=blue!30},
	aggreg/.style = {rectangle, minimum width=0.9cm, minimum height=0.75cm, text centered, text width=0.9cm, draw=black, fill=blue!30},
        aggregg/.style = {rectangle, minimum width=1cm, minimum height=0.75cm, text centered, text width=1cm, draw=black, fill=blue!30},
	decision/.style = {diamond, minimum width=0.8cm, minimum height=0.5cm, text centered, text width=0.8cm, draw=black, fill=green!30},
	start/.style = {rectangle, minimum width=1cm, minimum height=0.75cm, text centered, text width=1.75cm, draw=black, fill=red!30},
	arrow/.style = {thick,->,>=stealth},
	myfitd/.style={draw,dashed,blue!40, inner xsep=10pt, inner ysep=10pt, rounded corners=5pt, line width=1pt},
	myfit/.style={draw,dashed,blue, inner xsep=10pt, inner ysep=10pt, rounded corners=5pt, line width=1pt},
	myfits/.style={draw,black, inner xsep=12pt, inner ysep=15pt, rounded corners=10pt, line width=1pt},
	mytitles/.style={draw,black, fill=black!10, inner sep=5pt, right, xshift=20pt, line width=1pt},
	mytitless/.style={draw,densely dashed,blue!40, fill=blue!5, inner sep=2.5pt, right, xshift=10pt, line width=1pt},
	mytitle/.style={draw,densely dashed,blue, fill=blue!10, inner sep=2.5pt, right, xshift=10pt, line width=1pt}
}

\pgfplotsset{
	compat=1.11,
	legend image code/.code={
		\draw[mark repeat=2,mark phase=2]
		plot coordinates {
			(0cm,0cm)
			(0cm,0cm)        
			(0.4cm,0cm)         
		};
	}
}

\usetikzlibrary{arrows.meta, shapes.geometric}
\usepackage{cite}

\emergencystretch=\maxdimen
\usepackage[frozencache,cachedir=.]{minted}
\usepackage{array}
\newcommand{\directgolib}{\texttt{DIRECTGOLib v2.0}}

\newcommand{\dgo}{\texttt{DIRECTGO}}
\newcommand{\directgen}{\texttt{GENDIRECT}}

\newcommand{\matlab}{\texttt{MATLAB}}

\newcommand{\direct}{\texttt{DIRECT}}

\newcommand{\birmin}{\texttt{BIRMIN}}
\newcommand{\dirmin}{\texttt{DIRMIN}}

\newcommand{\mrdirect}{\texttt{MrDIRECT}}

\newcommand{\dtcg}{\texttt{1-DTC-GL}}

\newcommand{\halrect}{\texttt{HALRECT-IA}}

\newcommand{\indexsett}{\mathbb{I}}
\newcommand{\PreserveBackslash}[1]{\let\temp=\\#1\let\\=\temp}
\newcolumntype{C}[1]{>{\PreserveBackslash\centering}p{#1}}
\newcolumntype{R}[1]{>{\PreserveBackslash\raggedleft}p{#1}}
\newcolumntype{L}[1]{>{\PreserveBackslash\raggedright}p{#1}}
\newcommand{\algrule}[1][.2pt]{\par\vskip.5\baselineskip\hrule height #1\par\vskip.5\baselineskip}

\newcommand{\nosemic}{\renewcommand{\@endalgocfline}{\relax}}
\newcommand{\dosemic}{\renewcommand{\@endalgocfline}{\algocf@endline}}
\let\oldnl\nl
\newcommand{\nonl}{\renewcommand{\nl}{\let\nl\oldnl}}
\definecolor{onyx}{rgb}{0.06, 0.06, 0.06}
\definecolor{sandstorm}{rgb}{0.93, 0.84, 0.25}
\definecolor{princetonorange}{rgb}{1.0, 0.56, 0.0}
\definecolor{sienna}{rgb}{0.53, 0.18, 0.09}
\definecolor{psychedelicpurple}{rgb}{0.87, 0.0, 1.0}
\definecolor{bg}{rgb}{0.95,0.95,0.95}
\definecolor{ao}{rgb}{0.0, 0.5, 0.0}
\definecolor{arsenic}{rgb}{0.23, 0.27, 0.29}
\definecolor{armygreen}{rgb}{0.29, 0.33, 0.13}
\definecolor{antiquebrass}{rgb}{0.8, 0.58, 0.46}
\definecolor{DarkRed}{rgb}{0.55, 0.0, 0.0}
\definecolor{darkblue}{rgb}{0.0, 0.0, 0.55}
\definecolor{blueryb}{rgb}{0.01, 0.28, 1.0}
\definecolor{bluebell}{rgb}{0.64, 0.64, 0.82}
\definecolor{red}{rgb}{1.0, 0.0, 0.0}
\definecolor{redwood}{rgb}{0.67, 0.31, 0.32}
\definecolor{rose}{rgb}{1.0, 0.0, 0.5}
\definecolor{rosybrown}{rgb}{0.74, 0.56, 0.56}
\definecolor{rosewood}{rgb}{0.4, 0.0, 0.04}
\definecolor{saffron}{rgb}{0.96, 0.77, 0.19}
\definecolor{schoolbusyellow}{rgb}{1.0, 0.85, 0.0}
\definecolor{skyblue}{rgb}{0.53, 0.81, 0.92}
\definecolor{unmellowyellow}{rgb}{1.0, 1.0, 0.4}
\definecolor{wheat}{rgb}{0.96, 0.87, 0.7}
\definecolor{aureolin}{rgb}{0.99, 0.93, 0.0}
\definecolor{persianblue}{rgb}{0.11, 0.22, 0.73}
\definecolor{browna}{rgb}{0.59, 0.29, 0.0}
\definecolor{ufogreen}{rgb}{0.24, 0.82, 0.44}
\definecolor{forestgreen}{rgb}{0.13, 0.55, 0.13}
\definecolor{radicalred}{rgb}{1.0, 0.21, 0.37}
\definecolor{LightGreen}{rgb}{0.56, 0.93, 0.56}
\definecolor{LightCoral}{rgb}{0.94, 0.5, 0.5}
\definecolor{LightBlue}{rgb}{0.68, 0.85, 0.9}
\definecolor{DarkGreen}{rgb}{0.0, 0.2, 0.13}
\definecolor{LimeGreen}{rgb}{0.2, 0.8, 0.2}
\definecolor{DarkRed}{rgb}{0.55, 0.0, 0.0}
\definecolor{Tomato}{rgb}{1.0, 0.39, 0.28}
\definecolor{DarkBlue}{rgb}{0.0, 0.0, 0.55}
\definecolor{forestgreen}{rgb}{0.13, 0.55, 0.13}
\definecolor{zaffre}{rgb}{0.0, 0.08, 0.66}
\definecolor{wildstrawberry}{rgb}{1.0, 0.26, 0.64}
\definecolor{venetianred}{rgb}{0.78, 0.03, 0.08}
\definecolor{selectiveyellow}{rgb}{1.0, 0.73, 0.0}
\definecolor{yaleblue}{rgb}{0.06, 0.3, 0.57}

\pgfdeclareplotmark{aaa}{\node[star,star point ratio=2.25,minimum size=6pt,inner sep=0pt,draw=aureolin,solid, xshift=-0.5ex] {};}
\pgfdeclareplotmark{bbb}{\node[star,star points = 9,minimum size=6pt,inner sep=0pt,draw=DarkBlue,solid, xshift=-0.6ex] {};}
\pgfdeclareplotmark{ccc}{\node[isosceles triangle,draw,rotate=180,inner sep=0pt,draw=violet,solid,minimum size =5pt, xshift=-0.8ex] {};}
\pgfdeclareplotmark{ddd}{\node[isosceles triangle,draw,rotate=360,inner sep=0pt,draw=rosybrown,solid,minimum size =5pt, xshift=-0.5ex] {};}
\pgfdeclareplotmark{eee}{\node[isosceles triangle,draw,rotate=270,inner sep=0pt,draw=selectiveyellow,solid,minimum size =5pt, xshift=-1ex] {};}
\pgfdeclareplotmark{fff}{\node[star,star points = 7,minimum size=6pt,inner sep=0pt,draw=ufogreen,solid, xshift=-0.5ex] {};}
\pgfdeclareplotmark{ggg}{\node[regular polygon,regular polygon sides=5,minimum size=6pt,inner sep=0pt,draw=DarkRed,solid, xshift=-0.6ex] {};}
\pgfdeclareplotmark{hhh}{\node[regular polygon,regular polygon sides=7,minimum size=6pt,inner sep=0pt,draw=forestgreen,solid, xshift=-0.6ex] {};}
\pgfdeclareplotmark{iii}{\node[regular polygon,regular polygon sides=3,rotate=45,inner sep=0pt,draw=LightBlue,solid,minimum size =5pt, xshift=-0.4ex] {};}

\makeatletter
\tikzset{
	nomorepostactions/.code={\let\tikz@postactions=\pgfutil@empty},
	mymark/.style 2 args={decoration={markings,
			mark= between positions 0 and 1 step (1/11)*\pgfdecoratedpathlength with{%
				\tikzset{#2,every mark}\tikz@options
				\pgfuseplotmark{#1}%
			},
		},
		postaction={decorate},
		/pgfplots/legend image post style={
			mark=#1,mark options={#2},every path/.append style={nomorepostactions}
		},
	},
}
\makeatother

\definecolor{MutedBlue}{HTML}{1f77b4}
\definecolor{SafetyOrange}{HTML}{ff7f0e}
\definecolor{AsparagusGreen}{HTML}{2ca02c}
\definecolor{BrickRed}{HTML}{d62728}
\definecolor{MutedPurple}{HTML}{9467bd}
\definecolor{ChestnutBrown}{HTML}{8c564b}
\definecolor{RaspberryPink}{HTML}{e377c2}
\definecolor{MiddleGray}{HTML}{7f7f7f}
\definecolor{CurryYellowGreen}{HTML}{bcbd22}
\definecolor{BlueTeal}{HTML}{17becf}

\tikzstyle{s_dc}  = [mymark={o}{draw=black,solid,mark size=1.75pt},black,line width=0.75pt,opacity=0.8,thick]
\tikzstyle{s_dci} = [mymark={*}{draw=black,solid,mark size=1.75pt},black,line width=0.75pt,dashed,opacity=0.8,thick]
\tikzstyle{s_hl} = [mymark={square}{draw=blue,solid,mark size=1.75pt},blue,line width=0.75pt,opacity=0.8]
\tikzstyle{s_hli} = [mymark={square*}{draw=blue,solid,mark size=1.75pt},blue,line width=0.75pt,dashed,opacity=0.8]
\tikzstyle{s_mr} = [mymark={triangle}{draw=CurryYellowGreen,solid,mark size=1.75pt},CurryYellowGreen,line width=0.75pt,opacity=0.8,thick]
\tikzstyle{s_mri} = [mymark={triangle*}{draw=CurryYellowGreen,solid,mark size=1.75pt},CurryYellowGreen,line width=0.75pt,dashed,opacity=0.8,thick]
\tikzstyle{s_bm} = [mymark={diamond}{draw=red,solid,mark size=1.75pt},red,line width=0.75pt,opacity=0.8]
\tikzstyle{s_bmi} = [mymark={diamond*}{draw=red,solid,mark size=1.75pt},red,line width=0.75pt,dashed,opacity=0.8]
\tikzstyle{s_dm} = [mymark={x}{draw=RaspberryPink,solid,mark size=1.75pt},RaspberryPink,line width=0.75pt,opacity=0.8,thick]
\tikzstyle{s_dmi} = [mymark={+}{draw=RaspberryPink,solid,mark size=1.75pt},RaspberryPink,line width=0.75pt,dashed,opacity=0.8,thick]

\smartqed  

\begin{document}

\title{GENDIRECT: a GENeralized DIRECT-type algorithmic framework for derivative-free global optimization}

\titlerunning{GENDIRECT: A generalized DIRECT-type algorithmic framework for DFGO}


\author{Linas Stripinis \and Remigijus Paulavi\v{c}ius}

\institute{L. Stripinis \and R. Paulavi\v{c}ius \at
	Institute of Data Science and Digital Technologies, Vilnius University, Akademijos 4, LT-08663, Vilnius, Lithuania \\
	\email{linas.stripinis@mif.vu.lt}   \\
	R. Paulavi\v{c}ius \\
	\email{remigijus.paulavicius@mif.vu.lt}
}


\date{Received: date / Accepted: date}

\maketitle

\begin{abstract}

Over the past three decades, numerous articles have been published discussing the renowned \direct{} algorithm (DIvididing RECTangles).
These articles present innovative ideas to enhance its performance and adapt it to various types of optimization problems.
To consolidate and summarize this progress, we have recently introduced \dgo{}---a comprehensive collection featuring more than fifty deterministic, derivative-free algorithmic implementations based on the \direct{} framework.
\dgo{} empowers users to conveniently employ diverse \direct-type algorithms, enabling efficient solutions to practical optimization problems. 
Despite their variations, \direct-type algorithms share a common algorithmic structure and typically differ only at certain steps.

Recognizing this, we take further steps in generalization within this paper and propose \directgen---\texttt{GEN}eralized \direct-type framework that encompasses and unifies \direct-type algorithms under a single generalized approach.
\directgen{} offers a practical alternative to the creation of yet another ``new'' \direct-type algorithm that closely resembles existing ones. Instead, \directgen{} allows the efficient generation of known or novel \direct-type optimization algorithms by assembling different algorithmic components. 
This approach provides considerably more flexibility compared to both the \dgo{} toolbox and individual \direct-type algorithms.
In general, \directgen{} allows the creation of approximately a few hundred thousand combinations of \direct-type algorithms, facilitating user-friendly customization and the incorporation of new algorithmic components for further advancements.

By modifying specific components of five highly promising \direct-type algorithms found in the existing literature using \directgen{}, the significant potential of \directgen{} has been demonstrated. 
The resulting newly developed improved approaches exhibit greater efficiency and enhanced robustness in dealing with problems of varying complexity.

\end{abstract}

\section{Introduction}
\label{sec:Introduction}

Optimization problems encountered in scientific and engineering domains often involve objective functions that can only be obtained through ``black-box'' methods or simulations, lacking derivative information.
For example, Google's internal services frequently employ black-box optimization techniques with automated parameter tuning engines~\cite{Golovin2017}. 
Furthermore, objective function evaluations are becoming more computationally expensive as applications grow in size and complexity~\cite{Larson2019}. 
Consequently, calculating derivatives is often infeasible or impractical.
As a result, there is a growing emphasis on the development of derivative-free global optimization (DFGO) methods.
These methods are specifically designed to address the growing complexity and diversity of optimization problems, where derivative information is neither available nor practical to compute. 
This active development of DFGO methods addresses the need for efficient optimization techniques in scenarios where derivatives cannot be utilized.

This paper considers a box-constrained single-objective optimization problem
\begin{equation}\label{eq:opt-problem}
	\begin{aligned}
		& \min_{\mathbf{x}\in D} && f(\mathbf{x}),
	\end{aligned}
\end{equation}
where $f:\mathbb{R}^n \rightarrow \mathbb{R}$ is a potentially ``black-box'' Lipschitz-continuous objective function with an unknown Lipschitz constant, and $\mathbf{x} \in \mathbb{R}^n$ is the input vector of control variables.
Moreover, $f$ can be non-linear, multi-modal, non-convex, and non-differentiable.
We assume that $f$ can only be computed at any point of the feasible region, which is a $n$-dimensional hyper-rectangle
\[
D = [ \mathbf{a},  \mathbf{b}] = \{ \mathbf{x} \in \mathbb{R}^n: a_j \leq x_j \leq b_j, j = 1, \dots, n\}.
\]
However, there is no access to additional information on the objective function $f(\mathbf{x})$, such as gradients and the Hessian, as is typical for a ``black-box'' case. 

Among the solution techniques available for a given problem~\eqref{eq:opt-problem}, population-based meta-heuristic methods have gained widespread popularity.
Numerous approaches have been proposed and developed within this category~\cite{Agrawal2019}.
For global optimization problems that involve costly evaluations, model-based optimization algorithms are commonly employed.
Among these algorithms, Bayesian optimization~\cite{Jones1998} and various surrogate models\cite{Kudela2022} stand out as the leading state-of-the-art methods for optimizing expensive "black-box" functions.

\direct~\cite{Jones1993} presents an alternative specifically tailored for ``black-box'' global optimization by extending the classical Lipschitz optimization~\cite{Paulavicius2007,Paulavicius2010:ol,Pinter1996book,Sergeyev2011}, eliminating the requirement of knowing the Lipschitz constant.
In contrast to the stochastic methods discussed above, the \direct-type algorithms adhere to a deterministic pattern.
A recent comprehensive numerical benchmark study involving various derivative-free global optimization solvers\cite{Stripinis2022arx} highlighted that particularly for problems with lower dimensions, \direct-type algorithms can significantly outperform stochastic approaches. 
Furthermore, certain combinations of hybrid local search algorithms based on \direct-type methods and finite differences~\cite{Shi2021} demonstrated exceptional efficiency in solving high-dimensional problems. 
Consequently, designing and developing efficient \direct-type algorithms is crucial and driven by practical needs.

Inspired by these observations, we have recently introduced \dgo{}, a \matlab{} toolbox dedicated to DFGO. 
The latest release of \dgo{} includes a comprehensive collection of 52 distinct algorithmic implementations based on the \direct{} framework. 
However, recent empirical studies~\cite{Stripinis2022mdpi,Stripinis2022wcgo} have highlighted that even more efficient \direct-type algorithms can be achieved by innovatively combining existing algorithmic steps. 
It seems that many authors may not spend enough time exploring the most suitable algorithmic framework when developing and publishing new algorithms of type \direct.

Therefore, this study introduces a novel framework called \directgen{}, which offers a \texttt{GEN}eralized \direct-type approach to derivative-free global optimization.
\directgen{} enables the construction of any known or previously unexplored \direct-type algorithm. 
Instead of developing yet another ``new'' \direct-type algorithm, \directgen{} provides a rapid and effective way of combining different components to create customized \direct-type algorithms.

Using \directgen{}, users can identify and utilize the most suitable \direct-type algorithm for a given optimization problem based on the latest advances in the field. 
Compared to the \dgo{} toolbox and individual \direct-type algorithms, the \directgen{} framework offers a significantly higher level of flexibility. 
In fact, \directgen{} allows the design of a few hundred thousand combinations of \direct-type algorithms and facilitates user-friendly experimentation with new algorithmic components.

\directgen{} is implemented as a separate extension of \dgo{}, complemented by a dedicated graphical user interface (GUI). 
This GUI provides easy access to all the features and capabilities of \directgen{}, ensuring a seamless user experience.

The capability of \directgen{} is showcased by selecting five highly promising \direct-type algorithms from the existing literature, as identified in~\cite{Stripinis2022arx, Stripinis2022halrect}. 
By leveraging \directgen{}, specific components that were identified as weaknesses in these algorithms are modified. 
As a result, some of these algorithms demonstrate significantly improved efficiency, showcasing the potential of \directgen{} in optimizing and refining \direct-type algorithms using the most recent \directgolib.

This work makes several significant contributions, including
\begin{enumerate}
    \item Introduction of a novel framework called \directgen{}, which represents a \texttt{GEN}eralized \direct-type algorithmic framework.
    \item \directgen{} provides an efficient and innovative approach to generate \direct-type optimization algorithms, whether they are existing algorithms or entirely novel ones, by combining different algorithmic components.
    \item \directgen{} allows for the creation of a few hundred thousand combinations of \direct-type algorithms, facilitating user-friendly experimentation and enabling new developments in optimization.
    \item Description of the implementation of \directgen{} as an evolution of \dgo, complete with a separate graphical user interface (GUI) that ensures easy access to all its features. This implementation is free and open for anyone to use.
    \item Demonstration of the potential of \directgen{} by enhancing the efficiency of five chosen \direct-type algorithms through modifications. These modifications showcase the ability of \directgen{} to improve algorithmic performance further.
\end{enumerate}
In summary, this work contributes to the derivative-free global optimization field by introducing \directgen{}, a versatile framework that enables efficient algorithm generation, offers extensive customization options, and shows improved efficiency in established \direct-type algorithms.

The remainder of the paper is organized as follows.
\Cref{sec:rewiev} presents a concise overview of key advancements in the realm of \direct-type algorithms. 
\Cref{sec:directgen} introduces and elaborates on the \directgen{} framework. 
The experimental results of the newly developed algorithms and performance evaluation utilizing \directgen{} are analyzed in \Cref{sec:benchmarking}.
Lastly, \Cref{sec:conclusion} offers concluding remarks and outlines potential avenues for future exploration in this field.

\section{Background for \directgen{}}
\label{sec:rewiev}

\subsection{General structure of \direct-type algorithms}
\label{ssec:direct}

The \direct{} algorithm was originally designed to solve global optimization problems with box constraints~\eqref{eq:opt-problem}.
Despite numerous proposals, most follow a similar algorithmic structure and involve three primary steps: selection, sampling, and partitioning (see \Cref{alg:direct}).
However, at first, \direct-type algorithms typically transform a feasible region $D = [\mathbf{a}, \mathbf{b}]$ into a unit hyper-rectangle $ \bar{D} = [0, 1]^n$ referring to the original space $(D)$ solely to evaluate the objective function $f$ (as depicted in ~\Cref{alg:direct}, Lines~\ref{alg:initialization_begin}--\ref{alg:initialization_end}). 

The selection, partitioning, and sampling operations are executed within a normalized search domain $ \bar{D} $.
During each iteration, specific regions are identified as potentially optimal candidates (POC) and chosen for further investigation (see~\Cref{alg:direct}, Line~\ref{alg:selection_begin}).
In \direct-type algorithms, the objective function is sampled and evaluated at various points within each POC, which are then subdivided into smaller sub-regions (see Algorithm~\ref{alg:direct}, Lines~\ref{alg:sampling} and \ref{alg:subdivision}).
This selection, sampling, and subdivision process continues until a predefined limit is reached.

The subsequent subsections provide an overview of the primary techniques proposed for each step.
Although the selection step precedes sampling and partitioning, we will initially focus on the latter because the selection step relies directly on the strategies employed in sampling and partitioning.

\begin{algorithm}
    \normalsize
    \LinesNumbered
    \SetAlgoLined
    \SetKwInOut{Input}{input}
    \SetKwInOut{Output}{output}
    \SetKwData{Mmax}{M$_{\rm max}$}
    \SetKwData{Kmax}{K$_{\rm max}$}
    \SetKw{And}{and}
    \SetKw{Or}{or}
    
    \Input{Objective function $(f)$, search domain $(D)$, and adjustable algorithmic options $(opt)$: goal for the function value ($f^{\rm goal}$), maximal number of function evaluations ($\Mmax$) and algorithmic iterations ($\Kmax$) ; }
    \Output{The best found objective value $(f^{\min})$, solution point $(\mathbf{x}^{\min})$, and record of various performance metrics: percent error $(pe)$, number of iterations $(k)$, number of function evaluations $(m)$;} 
    \algrule
    
    \nonl \textbf{Initialization step:} \\ 
    \textit{Normalize} the search domain $D$ to the unit hyper-rectangle $\bar{D}$; \label{alg:initialization_begin} \\
    \textit{Evaluate} $f$ at initial sampling point(s) and set: \\
    $x^{\min}_j = \mid b_j - a_j \mid  \bar{x}_j + a_j, j=1, \dots, n$; \tcp*[f]{referring to $D$}\\
    $f^{\min} = f(\mathbf{x}^{\min})$; \\
    \textit{Initialize} performance measures: $k=1$, $m=1$ ;\label{alg:initialization_end}
    
    \While{$f^{\rm goal} < f^{\min}$ \And $m < \Mmax$ \And $k < \Kmax$ }{
        \textbf{Selection step:} \textit{Identify} the set $S_k$ of POCs\; \label{alg:selection_begin}
    
        \ForEach{$\bar{D}^j_k \in S_k$}{
            \textbf{Sampling step:} \textit{Evaluate} $f$ at newly sampled points in $\bar{D}^j_k$; \label{alg:sampling}\\
            \textbf{Partitioning step:} \textit{Subdivide} $\bar{D}^j_k$ \; \label{alg:subdivision}
        }\label{alg:global_end}
        \If{Hybrid}{
        \textit{Run} local search procedure;\tcp*[f]{only in hybrid version}\\
        }
        \textit{Update} $f^{\min}, \mathbf{x}^{\min}$, and performance measures: $k$ and $m$;
    }
    \textbf{Return} $f^{\min}, \mathbf{x}^{\min}$, and performance measures $(k, m)$.
    \caption{Main steps of \direct-type algorithms}
    \label{alg:direct}
\end{algorithm}

\subsection{Summary of sampling and partitioning schemes}
\label{sec:partitioning}

In this section, we present a brief summary of seven primary sampling and partitioning approaches that have been proposed in existing literature ~\cite{Jones1993,Jones2001,Sergeyev2006,Paulavicius2013:jogo,Paulavicius2016:jogo,Stripinis2022halrect} and implemented within the \directgen{} framework.
\Cref{tab:partitioning-strategies} provides an overview of these schemes, including illustrative examples from the initial iterations.
Blue-colored sub-regions indicate the POCs in the current partition.

Although each of the seven schemes possesses distinct characteristics, they demonstrate significant similarities. 
Particularly, these schemes involve sampling new points and subdividing larger regions into smaller, non-overlapping sub-regions. 
In cases where there is more than one longest side, two primary strategies for division emerge:
\begin{itemize}
    \item Subdivision along all dimensions with the maximum side length.
    \item Subdivision along a single dimension with maximum side length.
\end{itemize}
It is worth mentioning that the original \direct{} algorithm proposed subdividing along all dimensions. 
However, extensive experimentation has indicated that this approach does not consistently yield effective results.

\begin{table}[ht]
    \begin{minipage}{\textwidth}
        \caption{Summary of sampling and partitioning schemes commonly utilized in \direct-type algorithms implemented within \directgen{} (in ascending order of the year of publication)}
        \label{tab:partitioning-strategies}
        \begin{tabular}{p{0.129\textwidth}p{0.40\textwidth}p{0.37\textwidth}}
            \toprule
            \textbf{Notation \& Source} & Partitioning and sampling scheme & An example of the initialization and two subsequent iterations\\
            \midrule
            \textbf{DTC}~\cite{Jones2001}  & A hyper-rectangular partition based on one-\textbf{D}imensional \textbf{T}risection, and sampling points located at \textbf{C}enters. &
            \raisebox{-0.75\totalheight}{
            \begin{tikzpicture}
            \begin{groupplot}[group style={group size=3 by 1,horizontal sep=3pt}]            \nextgroupplot[width=0.25\textwidth,height=0.25\textwidth,ytick=\empty,xtick=\empty,enlargelimits=0]
            \draw [black, thick, mark size=0.1pt, fill=blue!30,line width=0.2mm] (axis cs:0, 0) rectangle (axis cs:1, 1);
            \addplot[thick,patch,mesh,draw,black,patch type=rectangle,line width=0.2mm] coordinates {(0, 0) (1, 0) (1, 1) (0, 1)} ;
            \addplot[only marks,mark=*, mark size=1.25pt,black] coordinates {(1/2, 1/2)};            \nextgroupplot[width=0.25\textwidth,height=0.25\textwidth,ytick=\empty,xtick=\empty,enlargelimits=0]
            \draw [black, thick, mark size=0.1pt, fill=blue!30,line width=0.2mm] (axis cs:2/3, 0) rectangle (axis cs:1, 1);
            \addplot[thick,patch,mesh,draw,black,patch type=rectangle,line width=0.2mm] coordinates {(0, 0) (1, 0) (1, 1) (0, 1)} ;
            \addplot[thick,patch,mesh,draw,black,patch type=rectangle,line width=0.2mm] coordinates {(1/3, 0) (1/3, 1) (2/3, 1) (2/3, 0)};
            \addplot[only marks,mark=o, mark size=1.25pt,black] coordinates {(1/2, 1/2)};
            \addplot[only marks,mark=o, mark size=1.25pt,black] coordinates {(1/6, 1/2)};
            \addplot[only marks,mark=*, mark size=1.25pt,blue]  coordinates {(5/6, 1/2)};            \nextgroupplot[width=0.25\textwidth,height=0.25\textwidth,ytick=\empty,xtick=\empty,enlargelimits=0]
            \draw [black, thick, mark size=0.1pt, fill=blue!30,line width=0.2mm] (axis cs:0, 0) rectangle (axis cs:1/3, 1);
            \draw [black, thick, mark size=0.1pt, fill=blue!30,line width=0.2mm] (axis cs:2/3, 1/3) rectangle (axis cs:1, 2/3);
            \addplot[thick,patch,mesh,draw,black,patch type=rectangle,line width=0.2mm] coordinates {(0, 0) (1, 0) (1, 1) (0, 1)} ;
            \addplot[thick,patch,mesh,draw,black,patch type=rectangle,line width=0.2mm] coordinates {(1/3, 0) (1/3, 1) (2/3, 1) (2/3, 0)};
            \addplot[thick,patch,mesh,draw,black,patch type=rectangle,line width=0.2mm] coordinates {(2/3, 1/3) (2/3, 2/3) (1, 2/3) (1, 1/3)};
            \addplot[only marks,mark=o, mark size=1.25pt,black] coordinates {(1/2, 1/2)};
            \addplot[only marks,mark=*, mark size=1.25pt,blue]  coordinates {(1/6, 1/2)};
            \addplot[only marks,mark=*, mark size=1.25pt,blue]  coordinates {(5/6, 1/2)};
            \addplot[only marks,mark=o, mark size=1.25pt,black] coordinates {(5/6, 1/6)};
            \addplot[only marks,mark=o, mark size=1.25pt,black] coordinates {(5/6, 5/6)};
            \end{groupplot}
            \end{tikzpicture}} \\
            \midrule
            \textbf{DTDV}~\cite{Sergeyev2006} & A hyper-rectangular partition based on one-\textbf{D}imensional \textbf{T}risection, and sampling points located at two \textbf{D}iagonal \textbf{V}ertices. &
            \raisebox{-0.75\totalheight}{
            \begin{tikzpicture}
            \begin{groupplot}[group style={group size=3 by 1,horizontal sep=2.5pt}]            \nextgroupplot[width=0.25\textwidth,height=0.25\textwidth,ytick=\empty,xtick=\empty,enlargelimits=0]
            \draw [black, thick, mark size=0.1pt, fill=blue!30,line width=0.2mm] (axis cs:0, 0) rectangle (axis cs:1, 1);
            \addplot[thick,patch,mesh,draw,black,patch type=rectangle,line width=0.2mm] coordinates {(0, 0) (1, 0) (1, 1) (0, 1)};
            \addplot[only marks,mark=*, mark size=1.25pt,blue]  coordinates {(0, 0)};
            \addplot[only marks,mark=*, mark size=1.25pt,blue] coordinates {(1, 1)};            \nextgroupplot[width=0.25\textwidth,height=0.25\textwidth,ytick=\empty,xtick=\empty,enlargelimits=0]
            \draw [black, thick, mark size=0.1pt, fill=blue!30,line width=0.2mm] (axis cs:0, 0) rectangle (axis cs:1/3, 1);
            \addplot[thick,patch,mesh,draw,black,patch type=rectangle,line width=0.2mm] coordinates {(0, 0) (1, 0) (1, 1) (0, 1)};
            \draw [black, thick, mark size=0.1pt,line width=0.2mm] (axis cs:1/3, 0) rectangle (axis cs:2/3, 1);
            \addplot[only marks,mark=*, mark size=1.25pt,blue]  coordinates {(0, 0)};
            \addplot[only marks,mark=*, mark size=1.25pt,blue]  coordinates {(1/3, 1)};
            \addplot[only marks,mark=o, mark size=1.25pt,black] coordinates {(2/3, 0)};
            \addplot[only marks,mark=o, mark size=1.25pt,black] coordinates {(1, 1)};            \nextgroupplot[width=0.25\textwidth,height=0.25\textwidth,ytick=\empty,xtick=\empty,enlargelimits=0]
            \draw [black, thick, mark size=0.1pt, fill=blue!30,line width=0.2mm] (axis cs:1/3, 0) rectangle (axis cs:2/3, 1);
            \draw [black, thick, mark size=0.1pt, fill=blue!30,line width=0.2mm] (axis cs:0, 0) rectangle (axis cs:1/3, 1/3);
            \addplot[thick,patch,mesh,draw,black,patch type=rectangle,line width=0.2mm] coordinates {(0, 0) (1, 0) (1, 1) (0, 1)} ;
            \draw [black, thick, mark size=0.1pt,line width=0.2mm] (axis cs:0, 1/3) rectangle (axis cs:1/3, 2/3);
            \draw [black, thick, mark size=0.1pt,line width=0.2mm] (axis cs:1/3, 0) rectangle (axis cs:2/3, 1);
            \addplot[only marks,mark=*, mark size=1.25pt,blue] coordinates {(0, 0)};
            \addplot[only marks,mark=*, mark size=1.25pt,blue] coordinates {(1/3, 1)};
            \addplot[only marks,mark=*, mark size=1.25pt,blue] coordinates {(2/3, 0)};
            \addplot[only marks,mark=o, mark size=1.25pt,black] coordinates {(1, 1)};
            \addplot[only marks,mark=o, mark size=1.25pt,black] coordinates {(0, 2/3)};
            \addplot[only marks,mark=*, mark size=1.25pt,blue] coordinates {(1/3, 1/3)};
            \end{groupplot}
            \end{tikzpicture}} \\
            \midrule
            \textbf{DTCS}~\cite{Paulavicius2013:jogo} & A simplicial partition based on one-\textbf{D}imensional \textbf{T}risection, and sampling points located at \textbf{C}enters of \textbf{S}implices. &
            \raisebox{-0.75\totalheight}{
            \begin{tikzpicture}
            \begin{groupplot}[group style={group size=3 by 1,horizontal sep=3pt}]            \nextgroupplot[width=0.25\textwidth,height=0.25\textwidth,ytick=\empty,xtick=\empty,enlargelimits=0]
            \addplot[thick,black,fill=blue!30,line width=0.2mm] coordinates {(0, 0) (0, 1) (1, 1)};
            \addplot[thick,patch,mesh,draw,black,patch type=rectangle,line width=0.2mm] coordinates {(0,0)(1,0)(1,1)(0,1)};
            \addplot[thick,patch,mesh,draw,black,line width=0.2mm] coordinates {(0, 0) (0, 1) (1, 1)};
            \addplot[only marks,mark=o, mark size=1.25pt,black]  coordinates {(2/3, 1/3)};
            \addplot[only marks,mark=*, mark size=1.25pt,blue]  coordinates {(1/3, 2/3)};            \nextgroupplot[width=0.25\textwidth,height=0.25\textwidth,ytick=\empty,xtick=\empty,enlargelimits=0]
            \addplot[thick,black,fill=blue!30,line width=0.2mm] coordinates {(0, 0) (1, 0) (1, 1)};
            \addplot[thick,black,fill=blue!30,line width=0.2mm] coordinates {(0, 0) (0, 1) (1/3, 1/3)};
            \addplot[thick,patch,mesh,draw,black,patch type=rectangle,line width=0.2mm] coordinates {(0,0)(1,0)(1,1)(0,1)};
            \addplot[thick,patch,mesh,draw,black,line width=0.2mm] coordinates {(1/3, 1/3) (0, 1) (2/3, 2/3)};
            \addplot[thick,patch,mesh,draw,black,line width=0.2mm] coordinates {(0, 0) (0, 1) (1, 1)};
            \addplot[only marks,mark=*, mark size=1.25pt,blue]  coordinates {(2/3, 1/3)};
            \addplot[only marks,mark=o, mark size=1.25pt,black]  coordinates {(1/3, 2/3)};
            \addplot[only marks,mark=*, mark size=1.25pt,blue]   coordinates {(1/9, 4/9)};
            \addplot[only marks,mark=o, mark size=1.25pt,black]  coordinates {(5/9, 8/9)}; \nextgroupplot[width=0.25\textwidth,height=0.25\textwidth,ytick=\empty,xtick=\empty,enlargelimits=0]
            \addplot[thick,black,fill=blue!30,line width=0.2mm] coordinates {(0, 0) (1/3, 1/3) (1, 0)};
            \addplot[thick,black,fill=blue!30,line width=0.2mm] coordinates {(0, 0) (1/3, 1/3) (0, 1/3)};
            \addplot[thick,patch,mesh,draw,black,patch type=rectangle,line width=0.2mm] coordinates {(0,0)(1,0)(1,1)(0,1)};
            \addplot[thick,patch,mesh,draw,black,line width=0.2mm] coordinates {(1/3, 1/3) (0, 1) (2/3, 2/3)};
            \addplot[thick,patch,mesh,draw,black,line width=0.2mm] coordinates {(1/3, 1/3) (1, 0) (2/3, 2/3)};
            \addplot[thick,patch,mesh,draw,black,line width=0.2mm] coordinates {(0, 0) (0, 1) (1, 1)};
            \addplot[thick,patch,mesh,draw,black,line width=0.2mm] coordinates {(0, 1/3) (1/3, 1/3) (0, 2/3)};
            \addplot[only marks,mark=o, mark size=1.25pt,black]  coordinates {(2/3, 1/3)};
            \addplot[only marks,mark=o, mark size=1.25pt,black]  coordinates {(1/3, 2/3)};
            \addplot[only marks,mark=*, mark size=1.25pt,blue]  coordinates {(4/9, 1/9)};
            \addplot[only marks,mark=o, mark size=1.25pt,black]  coordinates {(8/9, 5/9)};
            \addplot[only marks,mark=*, mark size=1.25pt,blue]  coordinates {(1/9, 2/9)};
            \addplot[only marks,mark=o, mark size=1.25pt,black]  coordinates {(1/9, 6/9)};
            \addplot[only marks,mark=o, mark size=1.25pt,black]  coordinates {(1/9, 4/9)};
            \addplot[only marks,mark=o, mark size=1.25pt,black]  coordinates {(5/9, 8/9)};
            \end{groupplot}
            \end{tikzpicture}} \\
            \midrule
            \textbf{DBVS}~\cite{Paulavicius2013:jogo} & A simplicial partition based on one-\textbf{D}imensional \textbf{B}isection, and sampling points located at \textbf{V}ertices of \textbf{S}implices. &
            \raisebox{-0.75\totalheight}{
            \begin{tikzpicture}
            \begin{groupplot}[group style={group size=3 by 1,horizontal sep=2.5pt}]            \nextgroupplot[width=0.25\textwidth,height=0.25\textwidth,ytick=\empty,xtick=\empty,enlargelimits=0]
            \addplot[thick,black,fill=blue!30,line width=0.2mm] coordinates {(0, 0) (0, 1) (1, 1)};
            \addplot[thick,patch,mesh,draw,black,patch type=rectangle,line width=0.2mm] coordinates {(0, 0) (1, 0) (1, 1) (0, 1)};
            \addplot[thick,patch,mesh,draw,black,line width=0.2mm] coordinates {(0, 0) (0, 1) (1, 1)};
            \addplot[only marks,mark=*, mark size=1.25pt,blue] coordinates {(0, 0)};
            \addplot[only marks,mark=*, mark size=1.25pt,blue] coordinates {(0, 1)};
            \addplot[only marks,mark=o, mark size=1.25pt,black] coordinates {(1, 0)};
            \addplot[only marks,mark=*, mark size=1.25pt,blue] coordinates {(1, 1)};            \nextgroupplot[width=0.25\textwidth,height=0.25\textwidth,ytick=\empty,xtick=\empty,enlargelimits=0]
            \addplot[thick,black,fill=blue!30,line width=0.2mm] coordinates {(0, 0) (1, 0) (1, 1)};
            \addplot[thick,black,fill=blue!30,line width=0.2mm] coordinates {(0, 0) (0, 1) (1/2, 1/2)};
            \addplot[thick,patch,mesh,draw,black,patch type=rectangle,line width=0.2mm] coordinates {(0, 0) (1, 0) (1, 1) (0, 1)};
            \addplot[thick,patch,mesh,draw,black,line width=0.2mm] coordinates {(0, 0) (0, 1) (1, 1)};
            \addplot[thick,patch,mesh,draw,black,line width=0.2mm] coordinates {(0, 0) (0, 1) (1/2, 1/2)};
            \addplot[only marks,mark=*, mark size=1.25pt,blue] coordinates {(0, 0)};
            \addplot[only marks,mark=*, mark size=1.25pt,blue] coordinates {(0, 1)};
            \addplot[only marks,mark=*, mark size=1.25pt,blue] coordinates {(1, 0)};
            \addplot[only marks,mark=*, mark size=1.25pt,blue] coordinates {(1, 1)};
            \addplot[only marks,mark=*, mark size=1.25pt,blue] coordinates {(1/2, 1/2)};            \nextgroupplot[width=0.25\textwidth,height=0.25\textwidth,ytick=\empty,xtick=\empty,enlargelimits=0]
            \addplot[thick,black,fill=blue!30,line width=0.2mm] coordinates {(0, 0) (1, 0) (1/2, 1/2)};
            \addplot[thick,black,fill=blue!30,line width=0.2mm] coordinates {(0, 0) (0, 1/2) (1/2, 1/2)};
            \addplot[thick,patch,mesh,draw,black,patch type=rectangle,line width=0.2mm] coordinates {(0, 0) (1, 0) (1, 1) (0, 1)};
            \addplot[thick,patch,mesh,draw,black,line width=0.2mm] coordinates {(0, 0) (0, 1) (1, 1)};
            \addplot[thick,patch,mesh,draw,black,line width=0.2mm] coordinates {(0, 0) (0, 1) (1, 0)};
            \addplot[thick,patch,mesh,draw,black,line width=0.2mm] coordinates {(0, 0) (0, 1/2) (1/2, 1/2)};
            \addplot[only marks,mark=*, mark size=1.25pt,blue] coordinates {(0, 0)};
            \addplot[only marks,mark=o, mark size=1.25pt,black] coordinates {(0, 1)};
            \addplot[only marks,mark=*, mark size=1.25pt,blue] coordinates {(1, 0)};
            \addplot[only marks,mark=o, mark size=1.25pt,black] coordinates {(1, 1)};
            \addplot[only marks,mark=*, mark size=1.25pt,blue] coordinates {(1/2, 1/2)};
            \addplot[only marks,mark=*, mark size=1.25pt,blue] coordinates {(0, 1/2)};
            \end{groupplot}
            \end{tikzpicture}} \\
            \midrule
            \textbf{DBDP}~\cite{Paulavicius2016:jogo} & A hyper-rectangular partition based on one-\textbf{D}imensional \textbf{B}isection, and sampling points located at two \textbf{D}iagonal \textbf{P}oints equidistant between themselves and a diagonal's vertices. &
            \raisebox{-0.9\totalheight}{
            \begin{tikzpicture}
            \begin{groupplot}[group style={group size=3 by 1,horizontal sep=3pt}]            \nextgroupplot[width=0.25\textwidth,height=0.25\textwidth,ytick=\empty,xtick=\empty,enlargelimits=0]
            \draw [black, thick,fill=blue!30, mark size=0.1pt,line width=0.2mm] (axis cs:0, 0) rectangle (axis cs:1, 1);
            \addplot[thick,patch,mesh,draw,black,patch type=rectangle,line width=0.2mm] coordinates {(0, 0) (1, 0) (1, 1) (0, 1)};
            \addplot[only marks,mark=*, mark size=1.25pt,blue] coordinates {(1/3, 1/3)};
            \addplot[only marks,mark=*, mark size=1.25pt,blue] coordinates {(2/3, 2/3)};            \nextgroupplot[width=0.25\textwidth,height=0.25\textwidth,ytick=\empty,xtick=\empty,enlargelimits=0]
            \draw [black, thick,fill=blue!30, mark size=0.1pt,line width=0.2mm] (axis cs:0, 0) rectangle (axis cs:1/2, 1);
            \addplot[thick,patch,mesh,draw,black,patch type=rectangle,line width=0.2mm] coordinates {(0, 0) (1, 0) (1, 1) (0, 1)};
            \draw [black, thick, mark size=0.1pt,line width=0.2mm] (axis cs:0, 0) rectangle (axis cs:1/2, 1);
            \addplot[only marks,mark=*, mark size=1.25pt,blue] coordinates {(1/3, 1/3)};
            \addplot[only marks,mark=o, mark size=1.25pt,black] coordinates {(2/3, 2/3)};
            \addplot[only marks,mark=o, mark size=1.25pt,black] coordinates {(5/6, 1/3)};
            \addplot[only marks,mark=*, mark size=1.25pt,blue] coordinates {(1/6, 2/3)};            \nextgroupplot[width=0.25\textwidth,height=0.25\textwidth,ytick=\empty,xtick=\empty,enlargelimits=0]
            \draw [black, thick,fill=blue!30, mark size=0.1pt,line width=0.2mm] (axis cs:1/2, 0) rectangle (axis cs:1, 1);
            \draw [black, thick,fill=blue!30, mark size=0.1pt,line width=0.2mm] (axis cs:0, 0) rectangle (axis cs:1/2, 1/2);
            \addplot[thick,patch,mesh,draw,black,patch type=rectangle,line width=0.2mm] coordinates {(0, 0) (1, 0) (1, 1) (0, 1)};
            \draw [black, thick, mark size=0.1pt,line width=0.2mm] (axis cs:0, 0) rectangle (axis cs:1/2, 1);
            \addplot[only marks,mark=*, mark size=1.25pt,blue] coordinates {(1/3, 1/3)};
            \addplot[only marks,mark=*, mark size=1.25pt,blue] coordinates {(2/3, 2/3)};
            \addplot[only marks,mark=*, mark size=1.25pt,blue] coordinates {(5/6, 1/3)};
            \addplot[only marks,mark=o, mark size=1.25pt,black] coordinates {(1/6, 2/3)};
            \addplot[only marks,mark=*, mark size=1.25pt,blue] coordinates {(1/6, 1/6)};
            \addplot[only marks,mark=o, mark size=1.25pt,black] coordinates {(1/3, 5/6)};
            \end{groupplot}
            \end{tikzpicture}} \\
            \midrule
            \textbf{DBVD}~\cite{Lakhdar2023} & A hyper-rectangular partition based on one-\textbf{D}imensional \textbf{B}isection, and sampling points located at one \textbf{V}ertice and one \textbf{D}iagonal point with a 2:3 diagonal distance from the sampling vertice. &
            \raisebox{-0.95\totalheight}{
            \begin{tikzpicture}
            \begin{groupplot}[group style={group size=3 by 1,horizontal sep=2.5pt}]            \nextgroupplot[width=0.25\textwidth,height=0.25\textwidth,ytick=\empty,xtick=\empty,enlargelimits=0]
            \draw [black, thick,fill=blue!30, mark size=0.1pt,line width=0.2mm] (axis cs:0, 0) rectangle (axis cs:1, 1);
            \addplot[thick,patch,mesh,draw,black,patch type=rectangle,line width=0.2mm] coordinates {(0, 0) (1, 0) (1, 1) (0, 1)};
            \addplot[only marks,mark=*, mark size=1.25pt,blue] coordinates {(1/3, 1/3)};
            \addplot[only marks,mark=*, mark size=1.25pt,blue] coordinates {(1, 1)};            \nextgroupplot[width=0.25\textwidth,height=0.25\textwidth,ytick=\empty,xtick=\empty,enlargelimits=0]
            \draw [black, thick,fill=blue!30, mark size=0.1pt,line width=0.2mm] (axis cs:0, 0) rectangle (axis cs:1/2, 1);
            \addplot[thick,patch,mesh,draw,black,patch type=rectangle,line width=0.2mm] coordinates {(0, 0) (1, 0) (1, 1) (0, 1)};
            \draw [black, thick, mark size=0.1pt,line width=0.2mm] (axis cs:0, 0) rectangle (axis cs:1/2, 1);
            \addplot[only marks,mark=*, mark size=1.25pt,blue] coordinates {(1/3, 1/3)};
            \addplot[only marks,mark=o, mark size=1.25pt,black] coordinates {(1, 1)};
            \addplot[only marks,mark=o, mark size=1.25pt,black] coordinates {(2/3, 1/3)};
            \addplot[only marks,mark=*, mark size=1.25pt,blue] coordinates {(0, 1)};            \nextgroupplot[width=0.25\textwidth,height=0.25\textwidth,ytick=\empty,xtick=\empty,enlargelimits=0]
            \draw [black, thick,fill=blue!30, mark size=0.1pt,line width=0.2mm] (axis cs:1/2, 0) rectangle (axis cs:1, 1);
            \draw [black, thick,fill=blue!30, mark size=0.1pt,line width=0.2mm] (axis cs:0, 0) rectangle (axis cs:1/2, 1/2);
            \addplot[thick,patch,mesh,draw,black,patch type=rectangle,line width=0.2mm] coordinates {(0, 0) (1, 0) (1, 1) (0, 1)};
            \draw [black, thick, mark size=0.1pt,line width=0.2mm] (axis cs:0, 0) rectangle (axis cs:1/2, 1);
            \addplot[only marks,mark=*, mark size=1.25pt,blue] coordinates {(1/3, 1/3)};
            \addplot[only marks,mark=*, mark size=1.25pt,blue] coordinates {(1, 1)};
            \addplot[only marks,mark=*, mark size=1.25pt,blue] coordinates {(2/3, 1/3)};
            \addplot[only marks,mark=o, mark size=1.25pt,black] coordinates {(0, 1)};
            \addplot[only marks,mark=*, mark size=1.25pt,blue] coordinates {(0, 0)};
            \addplot[only marks,mark=o, mark size=1.25pt,black] coordinates {(1/3, 2/3)};
            \end{groupplot}
            \end{tikzpicture}} \\
            \midrule
            \textbf{DBC}~\cite{Stripinis2022halrect} & A hyper-rectangular partition based on one-\textbf{D}imensional \textbf{B}isection, and sampling points located at \textbf{C}enter points. &
            \raisebox{-0.75\totalheight}{
            \begin{tikzpicture}
            \begin{groupplot}[group style={group size=3 by 1,horizontal sep=3pt}]            \nextgroupplot[width=0.25\textwidth,height=0.25\textwidth,ytick=\empty,xtick=\empty,enlargelimits=0]
            \draw [black, thick,fill=blue!30, mark size=0.1pt,line width=0.2mm] (axis cs:0, 0) rectangle (axis cs:1, 1);
            \addplot[thick,patch,mesh,draw,black,patch type=rectangle,line width=0.2mm] coordinates {(0, 0) (1, 0) (1, 1) (0, 1)};
            \addplot[only marks,mark=*, mark size=1.25pt,blue] coordinates {(1/2, 1/2)};            \nextgroupplot[width=0.25\textwidth,height=0.25\textwidth,ytick=\empty,xtick=\empty,enlargelimits=0]
            \draw [black, thick,fill=blue!30, mark size=0.1pt,line width=0.2mm] (axis cs:0, 0) rectangle (axis cs:1/2, 1);
            \addplot[thick,patch,mesh,draw,black,patch type=rectangle,line width=0.2mm] coordinates {(0, 0) (1, 0) (1, 1) (0, 1)};
            \draw [black, thick, mark size=0.1pt,line width=0.2mm] (axis cs:0, 0) rectangle (axis cs:1/2, 1);
            \addplot[only marks,mark=*, mark size=1.25pt,blue] coordinates {(1/2, 1/2)};
            \addplot[only marks,mark=*, mark size=1.25pt,blue] coordinates {(1/4, 1/2)};
            \addplot[only marks,mark=o, mark size=1.25pt,black] coordinates {(3/4, 1/2)};            \nextgroupplot[width=0.25\textwidth,height=0.25\textwidth,ytick=\empty,xtick=\empty,enlargelimits=0]
            \draw [black, thick,fill=blue!30, mark size=0.1pt,line width=0.2mm] (axis cs:1/2, 0) rectangle (axis cs:1, 1);
            \draw [black, thick,fill=blue!30, mark size=0.1pt,line width=0.2mm] (axis cs:0, 0) rectangle (axis cs:1/2, 1/2);
            \addplot[thick,patch,mesh,draw,black,patch type=rectangle,line width=0.2mm] coordinates {(0, 0) (1, 0) (1, 1) (0, 1)};
            \draw [black, thick, mark size=0.1pt,line width=0.2mm] (axis cs:0, 0) rectangle (axis cs:1/2, 1);
            \addplot[only marks,mark=*, mark size=1.25pt,blue] coordinates {(1/2, 1/2)};
            \addplot[only marks,mark=*, mark size=1.25pt,blue] coordinates {(1/4, 1/2)};
            \addplot[only marks,mark=*, mark size=1.25pt,blue] coordinates {(3/4, 1/2)};
            \addplot[only marks,mark=*, mark size=1.25pt,blue] coordinates {(1/4, 1/4)};
            \addplot[only marks,mark=o, mark size=1.25pt,black] coordinates {(1/4, 3/4)};
            \end{groupplot}
            \end{tikzpicture}} \\
            \bottomrule
        \end{tabular}
    \end{minipage}
\end{table}

\subsection{Selection schemes}
\label{sec:selection}

Initially, selecting POCs is straightforward since only one candidate is available, the entire feasible region.
However, to introduce selection schemes in subsequent iterations, we must first establish the concept of the current partition $(\mathcal{P}_k)$, which in iteration $k$, is defined as
\[
\mathcal{P}_k = \{ \bar{D}^i_k : i \in \indexsett_k \}, 
\]
where $ \bar{D}^i_k $ are hyper-rectangles (or simplices) and $ \indexsett_k $ is the index set identifying the current partition $ \mathcal{P}_k $.
Then, the next partition, $\mathcal{P}_{k+1}$, is obtained by subdividing the selected POCs from the current partition $ \mathcal{P}_k $.

When identifying POCs, two crucial aspects come into play: the $\bar{D}^i_k$ measure $(\bar{\delta}_k^i)$ and the general quality based on the function values attained at the sample points.

\subsubsection{Evaluating goodness of candidates}

The values of the objective function obtained from the sampled points $\mathbb{H}^i_k$ are utilized to assess the overall quality of the candidate.
We refer to this value as the aggregated function value $(\mathcal{F}_k^i)$, which represents the goodness of $ \bar{D}_k^i$. 
In summary, four strategies have been presented to evaluate $\mathcal{F}_k^i$ in the literature~\cite{Stripinis2022halrect} as defined in \Cref{{def:AggrFunctionValues}}.

\begin{definition}(Aggregated function values)
    \label{def:AggrFunctionValues}
    Let:
    \begin{itemize}
        \item $\bar{\delta}_k^i $ is a measure of $ \bar{D}_k^i$;
        \item $\mathbf{x}_{\rm m}^i$ is a midpoint of $ \bar{D}_k^i$;
        \item $\mathbf{x}^{\rm min}$ is a currently best found minimum point;
        \item $\mathbb{H}^i_k$ is a representative sampling index set of all sample points within $\bar{D}^i_k$;
        \item $\textrm{card}(\mathbb{H}^i_k)$ is the cardinality of a set $\mathbb{H}^i_k$.
    \end{itemize}
    Then:
    \begin{itemize}
        \item \textit{Midpoint value based aggregated function value}:
        \begin{equation}
            \label{eq:eq1}
            \mathcal{F}^i_k = f(\mathbf{x}_{\rm m}^i), 
        \end{equation}
        \item \textit{Minimum value based aggregated function value}:
        \begin{equation}
            \label{eq:eq2}
            \mathcal{F}^i_k = \min_{j \in \mathbb{H}^i_k} f(\mathbf{x}^j) 
        \end{equation}
        \item \textit{Mean value based aggregated function value}:
        \begin{equation}
            \label{eq:eq3}
            \mathcal{F}^i_k = \dfrac{1}{\textrm{card}( \mathbb{H}^i_k )} \sum_{j=1}^{\textrm{card}( \mathbb{H}^i_k )} f(\mathbf{x}^j)
        \end{equation}
        \item \textit{Midpoint and minimum values based aggregated function value}: 
        \begin{equation}
            \label{eq:eq4}
            \mathcal{F}^i_k = \frac{1}{2}\left( \min_{j \in \mathbb{H}^i_k} f(\mathbf{x}^j) +  f(\mathbf{x}_{\rm m}^i) \right)
        \end{equation}
    \end{itemize}  
\end{definition}

The use of the $\mathcal{F}_k^i$ evaluation strategy depends on the specific sampling strategy being utilized. 
For example, in the case of the DTC and DTCS schemes, the midpoint value-based $\mathcal{F}_k^i$ is adopted since sampling is performed solely at a single midpoint.
However, when there are multiple sampling points per candidate, alternative strategies have been shown to have a significant influence, as shown in previous work~\cite{Stripinis2022}.

\subsubsection{Measuring candidates}
\label{sec:measure}

Depending on the sampling strategy, basically only two different ways have been proposed to measure POCs:
\begin{eqnarray}
    \bar{\delta}_k^i &=& \lambda d_k^i, \label{eq:diagonal}\\
    \bar{\delta}_k^i &=& \max_{j, l \in \mathbb{H}^i_k} \| \mathbf{x}^{j} - \mathbf{x}^{l} \|_2.\label{eq:longside}
\end{eqnarray}
where $\lambda \in [0, 1]$ and $d_k^i$ represents the Euclidean length of the diagonal of $\bar{D}^i_k$ diagonal.
A couple of significant points should be emphasized in this regard.
First, instead of relying solely on the Euclidean norm, alternative norms (e.g., $\| \cdot \|_{\infty}$) have been observed to yield favorable results, as noted in the work~\cite{Gablonsky2001}.
Second, different partitioning schemes have employed various values for $\lambda$.
For example, some schemes use $\lambda = 1$, which corresponds to the full length of the diagonal, as seen in~\cite{Stripinis2022halrect}, while others adopt $\lambda = 2/3$, as demonstrated in~\cite{Paulavicius2016:jogo}.
However, since $\lambda$ applies uniformly to all partition elements $\bar{D}^i_k$ and serves solely to counterbalance the selection of POC, the choice of $\lambda$ does not affect the performance of algorithms of type \direct.

\subsubsection{Summary of POC selection schemes}

To address the identified limitations of \direct-type algorithms, various POC selection schemes have been proposed.
\Cref{def:selSchemes} defines the four most widely used and implemented selection schemes in \directgen{}, while a summary of them is given in \Cref{tab:partitioning-strategies}.

\begin{definition}(Selection schemes)
    \label{def:selSchemes}
    Let:
    \begin{itemize}
        \item $\mathcal{F}_k^i$ denotes the aggregated function value for $ \bar{D}_k^i$;
        \item $\bar{\delta}_k^i $ is a measure of $ \bar{D}_k^i$;
        \item $\mathbb{I}_k^i \subseteq \mathbb{I}_k$ represents a subset of indices that correspond to elements of $\mathcal{P}_k$ sharing the same measure ($\bar{\delta}_k^i$). Additionally, $\mathbb{I}_k^{\rm min}$ contains the indices of elements with the smallest measure, $\bar{\delta}_k^{\rm min}$, while $\mathbb{I}_k^{\rm max}$ --- with the largest. 
    \end{itemize}
    Then:
    \begin{itemize}
        \item \textit{Original selection}: A candidate $ \bar{D}_k^j, j \in \indexsett_k $ is said to be potentially optimal if there exists some rate-of-change (Lipschitz) constant $ \tilde{L} > 0$ such that
        \begin{equation}
            \label{eqn:potOptRect1}
            \mathcal{F}_k^j - \tilde{L}\delta_k^j \leq \mathcal{F}_k^i - \tilde{L}\delta_k^i, \quad \forall i \in \indexsett_k, 
        \end{equation}
        \item \textit{Aggresive selection}: For each $ \mathbb{I}_k^i \ (\min \le i \le \max) $ select $ \bar{D}^j_k, j \in \indexsett_k^i $ with the lowest $\mathcal{F}_k^j$, i.e.,
        \begin{equation}
            \label{eqn:potOptRectAggr}
            \mathcal{F}^j_k \leq \mathcal{F}^l_k, \quad \forall l \in \indexsett_k^i. 
        \end{equation}
        \item \textit{Pareto selection}: Select all candidates $\bar{D}_k^i, i \in \indexsett_k$ who are not dominated, which means that there is no other candidate $\bar{D}_k^j, j \in \indexsett_k$ that satisfies the condition:
        \begin{equation}
            \label{eq:pareto}
            (\delta_k^j \geq \delta_k^i \wedge \mathcal{F}_k^j < \mathcal{F}_k^i) \vee (\delta_k^j > \delta_k^i \wedge \mathcal{F}_k^j \leq \mathcal{F}_k^i).
        \end{equation}
        \item \textit{Reduced Pareto selection}: Select $ \bar{D}_k^i$ with the lowest $\mathcal{F}_k^i$ and $ \bar{D}_k^j$ with the most extensive measure $\delta_k^j$, breaking ties in favor of a lower value of the aggregate function.
    \end{itemize}
\end{definition}
In summary, aggressive selection aims to choose a comprehensive set of candidates, ensuring that at least one candidate is selected from each group with different diameters $(\delta_k^i)$ while prioritizing candidates with the lowest aggregated function value. 
Then, the number of candidates selected through Pareto-based criteria tends to exceed the original selection strategy. 
However, this approach, which emphasizes exploring candidates of intermediate sizes, can lead to slower convergence, particularly when dealing with less complex optimization problems.
Therefore, the primary motivation behind introducing a reduced set of Pareto-optimal candidates was to address this issue.

It is important to note that when multiple equally good POC exist with the same $\delta_k^i$ and $\mathcal{F}^i_k$, two distinct selection strategies are available:
\begin{itemize}
    \item Select all equally good POC;
    \item Select only one with the highest index number.
\end{itemize}
Furthermore, selection schemes can integrate additional conditions to enhance the balance between local and global directions.
The subsequent subsection provides a detailed examination of these conditions.

\begin{table}[htbp]
  \begin{minipage}{\textwidth}
    \caption{Summary of selection schemes implemented in \directgen{}}
    \label{tab:selection-schemes}
    \begin{tabular}{p{0.15\textwidth}p{0.3\textwidth}p{0.45\textwidth}}
      \toprule
      \textbf{Notation \& Source} & Description & Illustration of POCs selection (blue points) using \direct{} on the same sample\\
      \midrule
      \textbf{Original}~\cite{Jones1993} & \textit{\textbf{O}riginal selection strategy}. Selects POCs based on the lower Lipschitz bound estimates for all possible Lipschitz constant values. &
      \raisebox{-0.9\totalheight}{
      \begin{tikzpicture}
        \begin{axis}[
          width=0.465\textwidth,height=0.3\textwidth,every axis/.append style={font=\small},
          ymin=-0.05,ymax=0.4,xmin=0,xmax=0.6,ytick distance=0.1,xtick distance=0.15,
          ylabel style={yshift=-0.1cm},xlabel style={yshift=0.2cm},clip=false,enlargelimits=0.04,
          xticklabels={$0$,$0.00$,$0.15$,$0.30$,$0.45$,$0.60$},yticklabels={,$0.0$,$0.1$,$0.2$,$0.3$,$0.4$},
          xlabel={$\delta^i_k$},xlabel style={at={(1.075,-0.075)}},ylabel={$\mathcal{F}_k^i$},
          ],
          \addplot[only marks,mark=o,mark size=2pt,black!60,domain y=0:0.5] table[x=D,y=F] {ddr1.txt};
          \addplot[only marks,mark=*,mark size=2pt,blue] table[x=DO,y=FO] {ddr1.txt};
        \end{axis}
      \end{tikzpicture}} \\
      \midrule
      \textbf{Aggressive}~\cite{Baker2000} & \textit{Aggressive selection strategy}. Selects  at least one candidate from each group of different diameters. &
      \raisebox{-0.9\totalheight}{
      	\begin{tikzpicture}
      		\begin{axis}[
      			width=0.465\textwidth,height=0.3\textwidth,every axis/.append style={font=\small},
      			ymin=-0.05,ymax=0.4,xmin=0,xmax=0.6,ytick distance=0.1,xtick distance=0.15,
      			ylabel style={yshift=-0.1cm},xlabel style={yshift=0.2cm},clip=false,enlargelimits=0.04,
      			xticklabels={$0$,$0.00$,$0.15$,$0.30$,$0.45$,$0.60$},yticklabels={,$0.0$,$0.1$,$0.2$,$0.3$,$0.4$},
      			xlabel={$\delta^i_k$},xlabel style={at={(1.075,-0.075)}},ylabel={$\mathcal{F}_k^i$},
      			],
      			\addplot[only marks,mark=o,mark size=2pt,black!60,domain y=0:0.5] table[x=D,y=F] {ddr1.txt};
      			\addplot[only marks,mark=*,mark size=2pt,blue] table[x=DA,y=FA] {ddr1.txt};
      		\end{axis}
      \end{tikzpicture}} \\
      \midrule
      \textbf{Pareto}~\cite{mockus_2011} & \textit{Pareto selection strategy}. Selects all candidates that are non-dominated on size (the higher, the better) and aggregated function value (the lower, the better). &
      \raisebox{-0.9\totalheight}{
      	\begin{tikzpicture}
      		\begin{axis}[
      			width=0.465\textwidth,height=0.3\textwidth,every axis/.append style={font=\small},
      			ymin=-0.05,ymax=0.4,xmin=0,xmax=0.6,ytick distance=0.1,xtick distance=0.15,
      			ylabel style={yshift=-0.1cm},xlabel style={yshift=0.2cm},clip=false,enlargelimits=0.04,
      			xticklabels={$0$,$0.00$,$0.15$,$0.30$,$0.45$,$0.60$},yticklabels={,$0.0$,$0.1$,$0.2$,$0.3$,$0.4$},
      			xlabel={$\delta^i_k$},xlabel style={at={(1.075,-0.075)}},ylabel={$\mathcal{F}_k^i$},
      			],
      			\addplot[only marks,mark=o,mark size=2pt,black!60,domain y=0:0.5] table[x=D,y=F] {ddr1.txt};
      			\addplot[only marks,mark=*,mark size=2pt,blue] table[x=DP,y=FP] {ddr1.txt};
      		\end{axis}
      \end{tikzpicture}} \\
      \midrule
      \textbf{Reduced Pareto}~\cite{Mockus2013} & \textit{Reduced Pareto selection strategy}. Selects only two candidates, the first and the last point on the Pareto front. &
      \raisebox{-0.9\totalheight}{
      	\begin{tikzpicture}
      		\begin{axis}[
      			width=0.465\textwidth,height=0.3\textwidth,every axis/.append style={font=\small},
      			ymin=-0.05,ymax=0.4,xmin=0,xmax=0.6,ytick distance=0.1,xtick distance=0.15,
      			ylabel style={yshift=-0.1cm},xlabel style={yshift=0.2cm},clip=false,enlargelimits=0.04,
      			xticklabels={$0$,$0.00$,$0.15$,$0.30$,$0.45$,$0.60$},yticklabels={,$0.0$,$0.1$,$0.2$,$0.3$,$0.4$},
      			xlabel={$\delta^i_k$},xlabel style={at={(1.075,-0.075)}},ylabel={$\mathcal{F}_k^i$},
      			],
      			\addplot[only marks,mark=o,mark size=2pt,black!60,domain y=0:0.5] table[x=D,y=F] {ddr1.txt};
      			\addplot[only marks,mark=*,mark size=2pt,blue] table[x=DR,y=FR] {ddr1.txt};
      		\end{axis}
      \end{tikzpicture}} \\
      \bottomrule
    \end{tabular}
  \end{minipage}
\end{table}

\subsubsection{Additional approaches for improved local and global POC selection.}
\label{sec:tehniques}

\textbf{Excessive local refinement reduction techniques.}
To protect the algorithm against excessive refinement around current local minima $f^{\rm min}$, the authors in the \direct{} literature~\cite{Jones1993,Finkel2006,Liu2013} proposed incorporating one of the following conditions along with Eq. \eqref{eqn:potOptRect1} in the original selection scheme:
\begin{eqnarray}
    \mathcal{F}_k^j - \tilde{L}\delta_k^j &\leq& f^{\rm min} - \varepsilon|f^{\rm min}|,\label{eqn:min}\\
    \mathcal{F}_k^j - \tilde{L}\delta_k^j &\leq& f^{\rm min} - \varepsilon|f^{\rm min} - f^{\rm median}|,\label{eqn:median}\\
    \mathcal{F}_k^j - \tilde{L}\delta_k^j &\leq& f^{\rm min} - \varepsilon|f^{\rm min} - f^{\rm average}|\label{eqn:average}.
\end{eqnarray}
Therefore, the lower Lipschitz bound of the POC must be lower than the current minimum value $(f^{\rm min})$ to at least some extent.
The parameter $\varepsilon$ plays a crucial role in determining the adjustment of the lower Lipschitz bound.
In the study conducted by~\cite{Jones1993}, favorable results were achieved using values of $\varepsilon$ ranging from $10^{-3}$ to $10^{-7}$, and a default value of $\varepsilon = 10^{-4}$ is suggested.
To reduce the sensitivity of the objective function to additive scaling, subtraction of the median value $(f^{\rm median})$ or the average $(f^{\rm average})$ value (as shown in Eqs. \eqref{eqn:median} and \eqref{eqn:average}) was proposed.

\textbf{Restart technique for the $\varepsilon$ parameter}. In~\cite{Finkel2004aa}, an adaptive scheme is introduced for the parameter $\varepsilon$ to prevent wasteful function evaluations in minor regions $\bar{D}^i_k$ where negligible improvements are expected.
The restart technique begins with $\varepsilon = 0$ and is maintained until an improvement is observed.
However, if there is no improvement for five consecutive iterations, it suggests a potential stagnation at a local optimum.
To address this, the algorithm switches to $\varepsilon = 0.01$.
Within 50 iterations, the restart technique returns to $\varepsilon = 0$ if an improvement is found or no progress is made.
If another 50 iterations pass without improvement, this indicates a possible discovery of the global minimum, requiring further refinement.

\textbf{Multi-level candidate selection using different $\varepsilon$ values}. In~\cite{Liu2015,Liu2015b}, two alternative multi-level techniques are proposed for the candidate selection procedure, involving three different levels:
\begin{itemize}
    \item Level 2: The \direct-type algorithm is executed with the usual settings, employing $\varepsilon = 10^{-5}$.
    \item Level 1: The selection is limited to $90\%$ of $\bar{D}^k_i \in \mathcal{P}^k$, excluding $10\%$ of the candidates with the largest measure. In this level, $\varepsilon = 10^{-7}$ is used.
    \item Level 0: The selection is limited to $10\%$ of the candidates with the largest measure, disregarding those excluded at level 1. Here, $\varepsilon = 0$ is used.
\end{itemize}
Both strategies are recommended in the study cycle through these levels using a combination of the ``W-cycle'': $21011012$.
One of the methods~\cite{Liu2015} employs a fixed $\varepsilon = 10^{-4}$ value at all levels, while the other~\cite{Liu2015b} adheres to the rules mentioned above.

\textbf{Globally-biased selection}. In the works~\cite{Paulavicius2014:jogo,Paulavicius2019:eswa}, a two-phase approach with global bias was introduced.
The algorithm effectively determines the adequacy of exploring a local optimum by employing the globally biased scheme.
It terminates the local phase (referred to as the ``usual'' phase) to prevent wasteful function evaluations by excessive local refinement. 
Upon stopping the usual phase, the algorithm seamlessly transitions into a global phase, wherein the hyper-rectangles chosen for further exploration must meet a minimum size requirement.
This globally biased phase continues until a better minimum point is discovered or a maximum number of ``global iterations'' is reached. 
Subsequently, the algorithm reverts back to the usual phase.
The search process alternates between these two phases, namely, the usual phase and the globally-biased phase, until a specified stopping condition is fulfilled.

\textbf{Two-phase (Global-Local) selection}. In the work~\cite{Stripinis2018a}, a two-phase selection approach has been introduced. 
This approach expands the set of previously obtained POCs by incorporating additional candidates based on their proximity to the current best minimum point $\mathbf{x}^{\rm min}$. 
This expansion is performed by conducting a selection process using calculated distances (instead of aggregated function values) between the current best minimum point and all other candidates:
\begin{equation}
    \label{eq:eq5}
    \mathcal{F}^i_k = \| \mathbf{x}_{\rm m}^i - \mathbf{x}^{\rm min} \|_2.
\end{equation}
By including candidates that are closer to the current minimum point, this step facilitates faster and more extensive exploration around the current minimum point.

\subsection{Acceleration through hybridization techniques}
\label{ssec:hybrid}

To our knowledge, three hybridization strategies have been proposed for \direct-type algorithms~\cite{Jones2001,Liuzzi2010,Holmstrom2004,Paulavicius2019:eswa}.

The first strategy, originally suggested by the author of the original \direct{}~\cite{Jones2001}, was later refined and improved in a work~\cite{Paulavicius2019:eswa}.
The concept behind this strategy involves performing a local search only when the algorithm achieves an improvement in the best current solution value, denoted $f^{\rm min}$. 
The best current solution $f^{\rm min}$ can be updated using a local search method or a more suitable direct-type algorithm that enables faster local refinement.

The second strategy~\cite{Holmstrom2004} operates similarly to the first one. 
However, instead of performing a local search from a single starting point, this strategy employs a clustering algorithm to identify multiple appropriate starting points. 
The following steps are executed within this suggested method:

\begin{itemize}
    \item The \direct-type algorithm is run for a fixed number of function evaluations, typically set at $100n + 1$ as the default.
    \item The sampled points are analyzed using an adaptive clustering algorithm to determine the optimal number of clusters. Subsequently, a local search is performed from the best point within each cluster.
    \item Additionally, the \direct-type algorithm is run again.
    \item If the \direct-type algorithm improves $f^{\rm min}$, a final local search is performed from the best point.
\end{itemize}

In the third, aggressive strategy~\cite{Liuzzi2010}, initiate a local search from the midpoint of each POC. 
However, this approach has faced significant criticism for potentially generating excessive local searches, as many starting points may converge to the same local optimum.

\section{\directgen{} optimization software}
\label{sec:directgen}

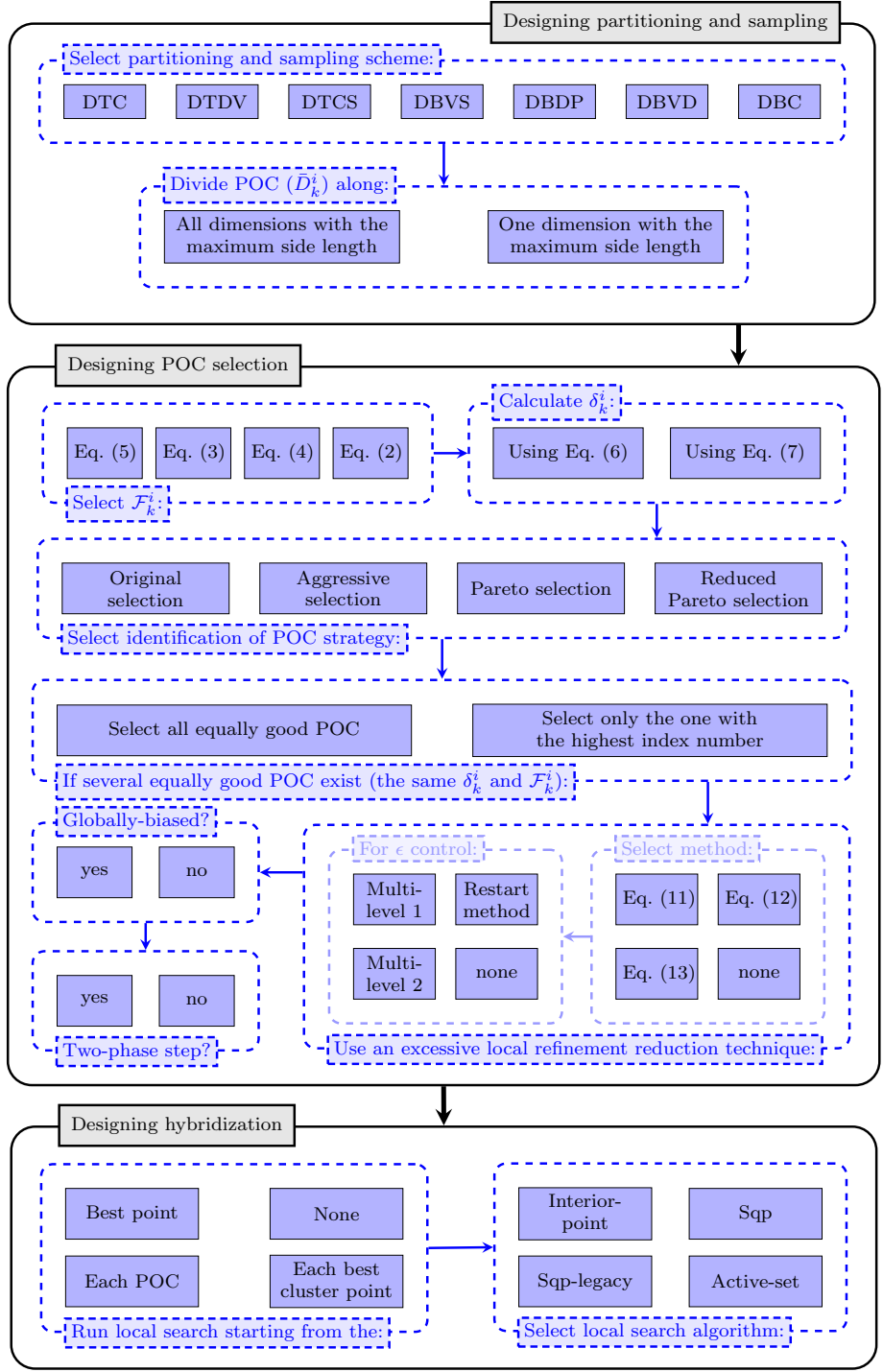
\begin{figure}[htbp]
	\centering
	\resizebox{0.99\textwidth}{!}{
		\begin{tikzpicture}[node distance=1.5cm]
			\node (p2) [scheme] {DTC};
			\node (p3) [scheme, right of=p2,xshift=0.15cm] {DTDV};
			\node (p4) [scheme, right of=p3,xshift=0.15cm] {DTCS};
			\node (p5) [scheme, right of=p4,xshift=0.15cm] {DBVS};
			\node (p6) [scheme, right of=p5,xshift=0.15cm] {DBDP};
			\node (p7) [scheme, right of=p6,xshift=0.15cm] {DBVD};
			\node (p8) [scheme, right of=p7,xshift=0.15cm] {DBC};
			\node[fit=(p2)(p3)(p4)(p5)(p6)(p7)(p8),myfit] (myfit1) {};
			\node[mytitle] at (myfit1.north west) {Select partitioning and sampling scheme:};
		
			\node (v1) [dimens, below of=p2,yshift=-0.5cm,xshift=2.6cm] {All dimensions with the maximum side length};
			\node (v2) [dimens, right of=v1,xshift=3.25cm] {One dimension with the maximum side length};
			\node[fit=(v1)(v2),myfit] (myfit2) {};
			\node[mytitle] at (myfit2.north west) {Divide POC $(\bar{D}^i_k)$ along:};
   
			\node[fit=(myfit1)(myfit2),myfits] (myfit66) {};
			\node[mytitles] at (myfit66.north) {Designing partitioning and sampling};

			\node (f1) [aggreg, below of=v1,xshift=1.3cm,yshift=-1.7cm] {Eq.~\eqref{eq:eq1}};
			\node (f2) [aggreg, left of=f1, xshift=0.2cm] {Eq.~\eqref{eq:eq3}};
			\node (f3) [aggreg, left of=f2, xshift=0.2cm] {Eq.~\eqref{eq:eq2}};
			\node (f4) [aggreg, left of=f3, xshift=0.2cm] {Eq.~\eqref{eq:eq4}};
			\node[fit=(f1)(f2)(f3)(f4),myfit] (myfit4) {};
			\node[mytitle] at (myfit4.south west) {Select $\mathcal{F}^i_k$:};

			\node (d1) [diamet, right of=f1,xshift=1.4cm] {Using Eq.~\eqref{eq:diagonal}};
			\node (d2) [diamet, right of=d1,xshift=1.1cm] {Using Eq.~\eqref{eq:longside}};
			\node[fit=(d1)(d2),myfit] (myfit3) {};
			\node[mytitle] at (myfit3.north west) {Calculate $\delta_k^i$:};

			\node (s1) [select, below of=f4,yshift=-0.5cm, xshift=0.6cm] {Original selection};
			\node (s2) [select, right of=s1,xshift=1.4cm] {Aggressive selection};
			\node (s3) [select, right of=s2,xshift=1.4cm] {Pareto selection};
			\node (s4) [select, right of=s3,xshift=1.4cm] {Reduced Pareto selection};
			\node[fit=(s1)(s2)(s3)(s4),myfit] (myfit6) {};
			\node[mytitle] at (myfit6.south west) {Select identification of POC strategy:};
			
			\node (e1) [selegh, below of=s1,yshift=-0.6cm,xshift=1.3cm] {Select all equally good POC};
			\node (e2) [selegh, right of=e1,xshift=4.6cm] {Select only the one with the highest index number};
			\node[fit=(e1)(e2),myfit] (myfit7) {};
			\node[mytitle] at (myfit7.south west) {If several equally good POC exist (the same $\delta_k^i$ and $\mathcal{F}^i_k$):};
			
			\node (r1) [aggregg, below of=e1,yshift=-1cm,xshift=6.2cm] {Eq.~\eqref{eqn:min}};
			\node (r2) [aggregg, right of=r1] {Eq.~\eqref{eqn:median}};
			\node (r3) [aggregg, below of=r1,yshift=0.4cm] {Eq.~\eqref{eqn:average}};
			\node (r4) [aggregg, right of=r3] {none};
			\node[fit=(r1)(r2)(r4)(r3),myfitd] (myfit8) {};
			\node[mytitless] at (myfit8.north west) {Select method:};
			
			\node (r5) [aggregg, left of=r1,xshift=-0.85cm] {Restart method};
			\node (r6) [aggregg, below of=r5,yshift=0.4cm] {none};
			\node (r7) [aggregg, left of=r5] {Multi-level $1$};
			\node (r8) [aggregg, left of=r6] {Multi-level $2$};
			\node[fit=(r5)(r6)(r7)(r8),myfitd] (myfit9) {};
			\node[mytitless] at (myfit9.north west) {For $\epsilon$ control:};
			
			\node[fit=(myfit8)(myfit9),myfit] (myfit10) {};
			\node[mytitle] at (myfit10.south west) {Use an excessive local refinement reduction technique:};
			
			\node (g1) [aggreg, left of=r7,xshift=-1.4cm,yshift=0.4cm] {no};
			\node (g2) [aggreg, left of=g1] {yes};
			\node[fit=(g1)(g2),myfit] (myfit11) {};
			\node[mytitle] at (myfit11.north west) {Globally-biased?};

			\node (l1) [aggreg, below of=g1, yshift=-0.4cm] {no};
			\node (l2) [aggreg, below of=g2, yshift=-0.4cm] {yes};
			\node[fit=(l1)(l2),myfit] (myfit5) {};
			\node[mytitle] at (myfit5.south west) {Two-phase step?};

			\node[fit=(myfit6)(myfit10)(myfit3)(myfit4),myfits] (myfit13) {};
			\node[mytitles] at (myfit13.north west) {Designing POC selection};

			\node (vv1) [diamest, below of=l1,yshift=-1.65cm,xshift=-0.95cm] {Best point};
			\node (vv3) [diamest, below of=vv1,yshift=0.5cm] {Each POC};
			\node (vv4) [diamest, right of=vv1,xshift=1.5cm] {None};
			\node (vv2) [diamest, below of=vv4,yshift=0.5cm] {Each best cluster point};
                \node[fit=(vv1)(vv2)(vv3)(vv4),myfit] (myfit222) {};
                \node[mytitle] at (myfit222.south west) {Run local search starting from the:};

                \node (vvv1) [diamest, right of=vv4,xshift=2.15cm] {Interior-point};
                \node (vvv2) [diamest, right of=vvv1,xshift=1cm] {Sqp};
                \node (vvv3) [diamest, right of=vv2,xshift=2.15cm] {Sqp-legacy};
                \node (vvv4) [diamest, right of=vvv3,xshift=1cm] {Active-set};
                \node[fit=(vvv1)(vvv2)(vvv3)(vvv4),myfit] (myfit2222) {};
                \node[mytitle] at (myfit2222.south west) {Select local search algorithm:};
                
                \node[fit=(myfit222)(myfit2222),myfits] (myfit21) {};
                \node[mytitles] at (myfit21.north west) {Designing hybridization};
				
			\draw[arrow,blue, line width=1pt] (myfit1.south-|myfit2.north) -- (myfit2.north);
			\draw[arrow,blue, line width=1pt] (myfit11.south) -- (myfit5.north);
			\draw[arrow,blue, line width=1pt] (myfit4.east) -- (myfit3.west);
			\draw[arrow,blue, line width=1pt] (myfit3.south) -- (myfit6.north-|myfit3.north);	
			\draw[arrow,blue, line width=1pt] (myfit6) -- (myfit7);
			\draw[arrow,blue, line width=1pt] (myfit7.south-|myfit8.north) -- (myfit10.north-|myfit8.north);
			\draw[arrow,blue!40, line width=1pt] (myfit8) -- (myfit9);
			\draw[arrow,blue, line width=1pt] (myfit11.west-|myfit10.west) -- (myfit11.east);	
			\draw[arrow,black, line width=2pt] (myfit66.south-|s4.north) -- (myfit13.north-|s4.north);
			\draw[arrow,black, line width=2pt] (myfit13.south) -- (myfit21.north);
            \draw[arrow,blue, line width=1pt] (myfit222.east) -- (myfit2222.west);
	\end{tikzpicture}}
	\caption{A flowchart for constructing \direct-type algorithm in \directgen.}
	\label{fig:directgens}
\end{figure}

This section describes the generalized algorithmic framework \direct{}. 
\Cref{fig:directgens}, illustrates the main architecture of the developed \directgen. 
Specifically, there are three large boxes in~\Cref{fig:directgens}, which represent the construction of the main \direct-type algorithmic steps within \directgen:
\begin{enumerate}
    \item The construction of partitioning and sampling scheme.
    \item The construction of the selection scheme.
    \item The construction of a hybridization scheme.
\end{enumerate}

The following subsections will provide a detailed exploration of how to effectively utilize \directgen{} using the \matlab{} command line interface and the dedicated graphical user interface (GUI).

\subsection{Utilizing \directgen{} through the command line interface.}
\label{ssec:command-line}

With \directgen{}, users can swiftly and effectively establish and solve global optimization problems by constructing a \direct-type algorithm via the \matlab{} command line interface. 
All relevant problem information is consolidated into a unified \matlab{} structure, which is then passed to the solver to extract the required data.

For the \directgen{} format, the solution process begins by generating the following structure:
\begin{minted}[bgcolor=bg]{matlab}
alg = GENDIRECT();
\end{minted}
The algorithm takes in a structured input that includes the optimization problem, dimension, lower and upper bounds, and a target value (if applicable).
Here is an example code snippet illustrating how these parameters can be set:
\begin{minted}[bgcolor=bg]{matlab}
alg.Problem.f     = 'objfun';    % Objective function
alg.Problem.n     = n;           % Dimension
alg.Problem.x_L   = zeros(n, 1); % Lower bounds
alg.Problem.x_U   = ones(n, 1);  % Upper bounds
alg.Problem.fgoal = 0.01;        % Optimal value set as target
alg.Problem.info  = false;       % Extract info from problem
\end{minted}
If the \texttt{alg.Problem.info} parameter is set to `true', the algorithm retrieves all the relevant information about the objective function from the `objfun' problem.

As we utilize test problems provided by \directgolib{}~\cite{DIRECTGOLib2022}, the stored information encompasses both the problem structure and the objective function. 
Consequently, the algorithms automatically extract all essential details from the given problem, including:
\begin{itemize}
    \item The dimensionality of the problem;
    \item The lower and upper bounds for each variable;
    \item The objective function value of the known solution;
    \item The solution point.
\end{itemize}
For further guidance on the utilization of \directgolib{}, additional information can be found in references~\cite{DIRECTGO2022,Stripinis2022}.

Users who want to customize the default algorithmic settings should utilize the \texttt{optParam} structure:
\begin{minted}[bgcolor=bg]{matlab}
alg.optParam.maxevals = 100;  % Maximal number of evaluations
alg.optParam.maxits   = 100;  % Maximal number of iterations
alg.optParam.showits  = true; % Show iteration status
\end{minted}

The next step involves constructing the algorithm using the procedures described in \Cref{tab:parameters}.

\begin{table}[ht]
    \caption{The parameters of \directgen{} used to construct \direct-type algorithms, with default values highlighted in \textcolor{blue}{blue}.}
    \label{tab:parameters}
    \begin{tabular}{llp{0.67\textwidth}}
        \toprule
        Step & Parameter & Description \\
        \midrule
        \multirow{4}{*}{\rotatebox{90}{Partitioning}} 
        & \texttt{Strategy} & Specify partitioning and sampling scheme (see \Cref{tab:partitioning-strategies}): \textcolor{blue}{\texttt{DTC}}, \texttt{DTDV}, \texttt{DTCS}, \texttt{DBVS}, \texttt{DBDP}, \texttt{DBVD}, or \texttt{DBC}.  \\
        \cmidrule{2-3}
        & \texttt{SubSides} & Specify subdivision strategy for multiple longest sides (see \Cref{sec:partitioning}):  \texttt{One} or \textcolor{blue}{\texttt{ALL}}.\\
        \midrule
        & \texttt{AggrFuncVal} & Specify strategy for a aggregated function value: \textcolor{blue}{\texttt{Midpoint}} (Eq.~\eqref{eq:eq1}), \texttt{Minimum} (Eq.~\eqref{eq:eq2}), \texttt{Mean} (Eq.~\eqref{eq:eq3}) or \texttt{MidMin} (Eq.~\eqref{eq:eq4}).\\
        \cmidrule{2-3}
        & \texttt{CandMeasure} & Specify strategy for a measure: \textcolor{blue}{\texttt{Diagonal}} (Eq.~\eqref{eq:diagonal}), or \texttt{LongSide} (Eq.~\eqref{eq:longside}).\\
        \cmidrule{2-3}
        & \texttt{Strategy} & Specify selection scheme (see \Cref{tab:selection-schemes}): \textcolor{blue}{\texttt{Original}}, \texttt{Aggressive}, \texttt{Pareto}, or \texttt{RedPareto}.  \\
        \cmidrule{2-3}
        \multirow{4}{*}{\rotatebox{90}{Selection}} & \texttt{EqualCand} & Specify behavior for equally good POC: \textcolor{blue}{\texttt{All}} or \texttt{One}.\\
        \cmidrule{2-3}
        & \texttt{SolRefin} & Specify excessive local refinement reduction technique: \textcolor{blue}{\texttt{Min}} (Eq.~\eqref{eqn:min}), \texttt{Median} (Eq.~\eqref{eqn:median}), \texttt{Average} (Eq.~\eqref{eqn:average}) or \texttt{Off}.\\
        \cmidrule{2-3}
        & \texttt{Ep} & Specify the value for $\varepsilon$ (Eqs.~\eqref{eqn:min}, \eqref{eqn:median}, \eqref{eqn:average}): \textcolor{blue}{$10^{-4}$}.\\
        \cmidrule{2-3}
        & \texttt{ControlEp} & Specify control technique for $\varepsilon$ (see \Cref{sec:tehniques}): \textcolor{blue}{\texttt{Off}}, \texttt{Restart}, \texttt{MultiLevel1} or \texttt{MultiLevel2}.\\
        \cmidrule{2-3}
        & \texttt{GloballyBiased} & Enable globally-biased POC selection (see \Cref{sec:tehniques}): \textcolor{blue}{\texttt{Off}} or \texttt{On}.\\
        \cmidrule{2-3}
        & \texttt{TwoPhase} & Enable two-phase selection of POC using \texttt{Distances} (Eq.~\eqref{eq:eq5})(see \Cref{sec:tehniques}): \textcolor{blue}{\texttt{Off}} or \texttt{On}.\\
        \midrule
        \multirow{10}{*}{\rotatebox{90}{Hybridization}}
        & \texttt{Strategy} & Specify hybridization strategy (see \Cref{ssec:hybrid}): \textcolor{blue}{\texttt{Off}}, \texttt{Single}, \texttt{Clustering} or \texttt{Aggressive}. \\
        \cmidrule{2-3}
        & \texttt{LocalSearch} & Specify derivative-free local search subroutine: \textcolor{blue}{\texttt{interior-point}}, \texttt{sqp}, \texttt{sqp-legacy} or \texttt{active-set}.\\
        \cmidrule{2-3}
        & \texttt{MaxIterations} & Specify the maximum iteration limit for a single local search subroutine call: \textcolor{blue}{$1000$}.\\
        \cmidrule{2-3}
        & \texttt{MaxEvaluations} & Specify the maximum function evaluation limit for a single local search subroutine call: \textcolor{blue}{$3000$}.\\
        \bottomrule
    \end{tabular}
\end{table}
After completing these steps, the algorithm is ready to solve the given problem using the following line of code:
\begin{minted}[bgcolor=bg]{matlab}
Results = alg.solve;
\end{minted}
Once the algorithm completes its computations, it returns the \texttt{Results} structure, which contains the optimization outcomes.

The subsequent subsections will outline the process of constructing \direct-type algorithmic steps.

\subsubsection{Designing partitioning and sampling scheme}
\label{sec:directgencombination}

To create a combination of \direct-type algorithms, the user needs to integrate components that determine the division and sampling strategy of the optimization domain.
The core framework for constructing the partitioning and sampling strategy is illustrated in the top block of \Cref{fig:directgens}.
The subsequent command lines illustrate how to configure the partitioning strategy of the original \direct{} algorithm:
\begin{minted}[bgcolor=bg]{matlab}
alg.Partitioning.Strategy = 'DTC';      
alg.Partitioning.SubSides = 'All'; 
\end{minted}
As a result of the given partitioning and sampling scheme in \Cref{fig:directgens}, there are $14$ possible combinations in \directgen.

\subsubsection{Designing the selection scheme}\label{sec:directgenselection}

Once the partitioning and sampling strategy has been established, the subsequent task is to determine the POC selection scheme.
Here is an example that illustrates the parameter values required for performing POC selection introduced in the original \direct{} algorithm:
\begin{minted}[bgcolor=bg]{matlab}
alg.Selection.AggrFuncVal = 'Midpoint';  
alg.Selection.CandMeasure = 'Diagonal'; 
alg.Selection.Strategy = 'Original';    
alg.Selection.EqualCand = 'All'; 
alg.Selection.SolRefin = 'Min';     
alg.Selection.Ep = 0.0001;    
alg.Selection.ControlEp = 'Off';   
alg.Selection.GloballyBased = 'Off';
alg.Selection.TwoPhase = 'Off';
\end{minted}

When the two-phase selection step is enabled, as demonstrated in the following code snippet:
\begin{minted}[bgcolor=bg]{matlab}
alg.Selection.TwoPhase = 'On';  
\end{minted}
the algorithm uses the designed selection scheme ('alg.Selection') to expand the set of promising candidate solutions (POC) based on the calculated distances obtained using Eq.~\eqref{eq:eq5}.
It is easy to calculate in \Cref{fig:directgens}, there are $4096$ different combinations for the selection steps of POC in \directgen.

\subsubsection{Designing hybridization scheme}
\label{sec:hybridirect}

In the third block of \Cref{fig:directgens}, users are required to select the desired hybridization technique. 
There are only $13$ possible combinations available in this block.

For example, to specify a hybridization scheme that utilizes a strategy calling an SQP local search (parameter \texttt{sqp}) subroutine only when an improvement in the best current solution is achieved (parameter \texttt{Single}), the following code can be used:
\begin{minted}[bgcolor=bg]{matlab}
alg.Hybridization.Strategy = 'Single';
alg.Hybridization.LocalSearch = 'sqp';
\end{minted}

\subsection{Utilizing \directgen{} through the graphical user interface}
\label{sec:gui}

\directgen{} is also accessible through the graphical user interface (GUI) of \dgo{}.
This GUI enables users to use \directgen{} without requiring prior programming or algorithmic knowledge. 
To access the \directgen{} tool, users can navigate to the \matlab{} APPS menu on the toolbar. 
Within \dgo{}, the generalized \direct{} algorithm (\directgen) can be selected from the algorithm drop-down menu.

\begin{figure}[ht]
    \centering
    \includegraphics[width=0.95\linewidth]{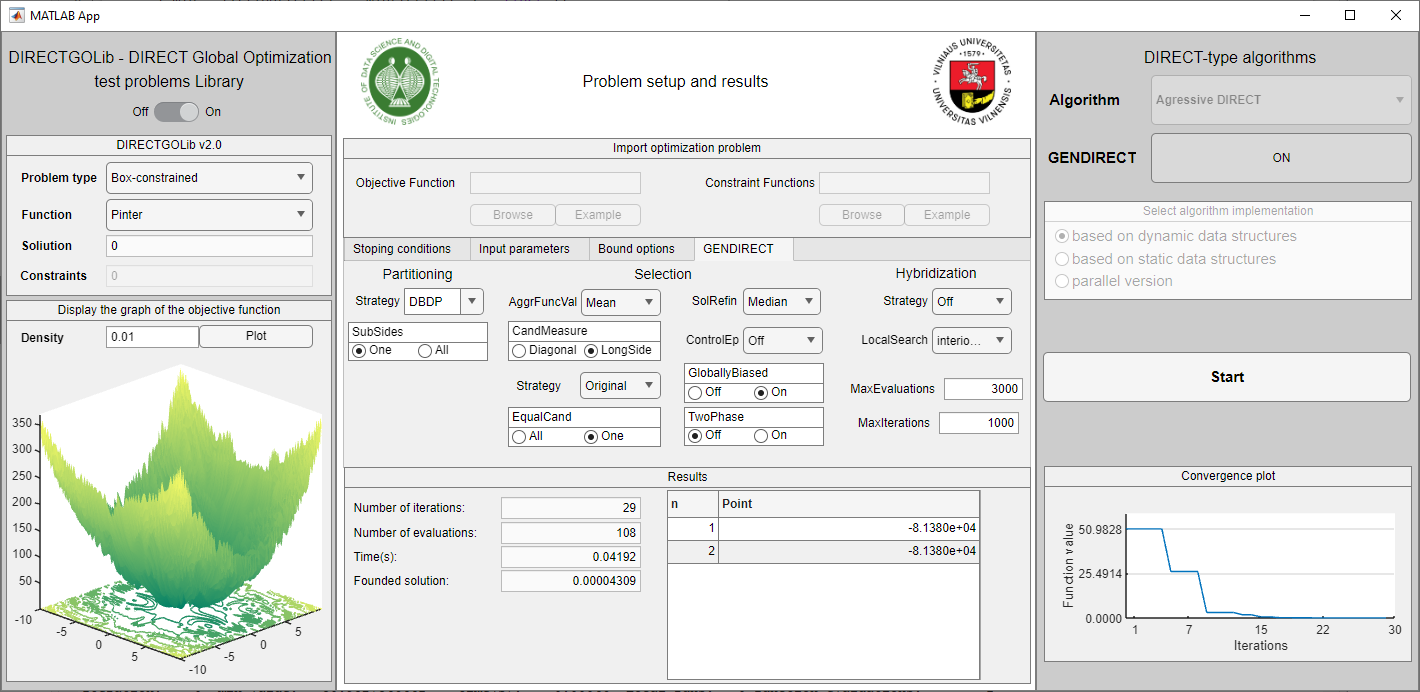}
    \caption{A snapshot of the graphical user interface (GUI) of \directgen{} in the \dgo{} software package.}
    \label{fig:directgo}
\end{figure}

The graphical interface of the main toolbox window \dgo{} is depicted in \Cref{fig:directgo}.
The \directgen{} window is centrally located in the GUI and facilitates the construction of the \direct{} algorithm by providing user-friendly functionalities.
For more comprehensive details of \dgo{}, see~\cite{Stripinis2022}.

\subsection{Remarks regarding the extension of \directgen{}}
\label{sec:expension}

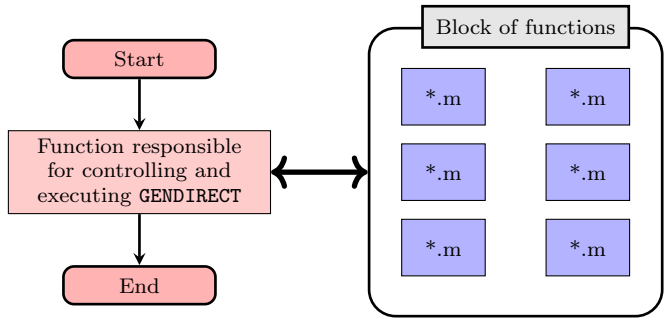
\begin{figure}[ht]
    \centering
    \begin{tikzpicture}[node distance=1cm]
        \node (p1) [startstop] {Start};
        \node (p2) [aggreg, right of=p1,xshift=3cm,yshift=-0.5cm] {*.m};
        \node (p4) [aggreg, below of=p2] {*.m};
        \node (p5) [aggreg, below of=p4] {*.m};
        \node (p12) [aggreg, right of=p2,xshift=0.9cm] {*.m};
        \node (p14) [aggreg, right of=p4,xshift=0.9cm] {*.m};
        \node (p15) [aggreg, right of=p5,xshift=0.9cm] {*.m};
        \node[fit=(p2)(p4)(p5)(p12),myfits] (myfit1) {};
        \node[mytitles] at (myfit1.north west) {Block of functions};
        \node (c1) [dimensss, below of=p1,yshift=-0.5cm] {Function responsible for controlling and executing \directgen};
        \node (c2) [startstop, below of=c1,yshift=-0.5cm] {End};
        
        \draw[arrow,black, line width=1pt] (p1) -- (c1);
        \draw[arrow,black, line width=1pt] (c1) -- (c2);
        \draw[<->,black, line width=2pt] (myfit1.west) -- (c1);
    \end{tikzpicture}
    \caption{The framework of the generalized \direct{} algorithm system (\directgen)}
    \label{fig:directrenas}
\end{figure}

\directgen{} comprises two primary components, as illustrated in \Cref{fig:directrenas}. 
Firstly, a function block encompasses various implementations of the steps involved in \direct-type algorithms. 
Secondly, the control structure ensures the seamless connection of algorithm components, facilitating the execution of the algorithm.

If a researcher intends to integrate a newly proposed step into \directgen, the function should be added to the function block. 
Ensuring that the implemented function adheres to the existing code's style is important. 
Subsequently, in the control function of \directgen, the newly created function should be incorporated accordingly, allowing \directgen{} to utilize it effectively.

\section{Simulation results and in-depth analysis}
\label{sec:benchmarking}
This section presents an analysis of the experimental results for newly developed improved algorithms and their performance evaluation using \directgen{}.

\subsection{An overview of benchmark test problems}

We employed a comprehensive set of $324$ benchmark test functions to thoroughly evaluate the newly proposed \directgen{} algorithm. 
These test problems were sourced from the latest version of the \directgolib{} library~\cite{Directgolib_2023}, which is built within the \matlab{} environment. 
The \directgolib{} integrates ten libraries and collections of well-established and recently developed test problems.

In \Cref{tab:directlibv2}, we present a summary of \directgolib{} and its constituent libraries. 
The table provides essential details, including references, publication years, the pool of problems, and the counts of scalable, separable, and multi-modal problems. 
Specifically, the table comprises $136$ test problems with fixed dimensions and $188$ test benchmarks that can be adjusted to any dimension size ($n$). 
For these test problems, we consider instances with variables set at $n = 2, 5,$ and $10$. 
However, it is worth noting that some functions, such as certain CEC functions~\cite{liang2013problem, wu2017problem}, are not applicable in all dimensions.

In our study, we thoroughly examined a total of $634$ test problems available in \directgolib{} to ensure comprehensive and robust evaluations of the proposed algorithmic framework \directgen{}.

\begin{table}[ht]
    \caption{Compilation of test problems from various libraries in the latest version of the \directgolib{} for box-constrained global optimization.}
    \begin{tabular*}{\textwidth}{@{\extracolsep{\fill}}lccccc}
        \toprule
        \multirow{2}{*}{Source} & \multirow{2}{*}{Year} & \multicolumn{4}{c}{Problems} \\
        \cmidrule{3-6}
        && Total & Scalable & Seperable & Multi-modal \\
        \midrule
        Hedar,~\cite{Hedar2005} & $2005$ & $31$ & $17$ & $8$ & $23$ \\
        Hansen~et~al.,\cite{Hansen2009bbob} & $2009$ & $24$ & $24$ & $5$ & $14$ \\
        Jamil~et~al.,\cite{Jamil2013} & $2013$ & $167$ & $69$ & $49$ & $127$ \\
        Gavana,~\cite{Gavana2021} & $2013$ & $193$ & $76$ & $64$ & $156$ \\
        Surjanovic~et~al.,~\cite{Surjanovic2013} & $2013$ & $50$ & $23$ & $11$ & $40$ \\
        Liang~et~al.,~\cite{liang2013problem} & $2014$ & $27$ & $27$ & $5$ & $24$ \\
        Wu~et~al.,~\cite{wu2017problem} & $2017$ & $20$ & $20$ & $0$ & $19$ \\
        Oldenhuis,~\cite{Oldenhuis2020} & $2020$ & $41$ & $12$ & $5$ & $33$ \\
        Layeb,~\cite{Layeb2022} & $2022$ & $18$ & $18$ & $2$ & $16$ \\
        Kudela~et~al.,~\cite{Kudela2022abs} & $2022$ & $8$ & $8$ & $8$ & $8$ \\
        \midrule
        Stripinis~et~al.,~\cite{Directgolib_2023} & $2023$ & $324$ & $188$ & $97$ & $261$ \\
        \bottomrule
    \end{tabular*}
    \label{tab:directlibv2}
\end{table}

In order to ensure that the global minimum point does not coincide with the initial sampling point in any tested algorithm, we employ shift operations. 
In other words, we randomly shift the solutions in the $X$-space. 
This involves transforming a given point $\mathbf{x}$ into $\hat{\mathbf{x}}$ using the following equation:
\begin{equation}
    \label{eq:shift}
    \hat{x}_j = \min \left\{ \max \left\{ x_j - \rho_j \lambda \vec{x}_j, a_j \right\}, b_j \right\},~j = 1,...,n.
\end{equation}
Here, $\vec{\mathbf{x}}$ is a randomly distributed random direction vector generated using the Mersenne-Twister pseudorandom generator, and $\lambda$ is a step size that serves two important purposes:
\begin{itemize}
    \item It prevents the global optima from moving outside of the feasible region.
    \item It allows for a more efficient placement of the solution within the problem domain, considering that different problems may have significantly different domain sizes.
\end{itemize}
The value of $\lambda$ is calculated by solving the following linear programming problem:
\begin{equation}
\label{eq:shiftlp}
    \begin{aligned}
        \max \,\,    & \lambda \\ 
        \text{s.t. } & \mathbf{x}^* + \lambda \vec{\mathbf{x}} \geq \mathbf{a}\\
                     & \mathbf{x}^* + \lambda \vec{\mathbf{x}} \leq \mathbf{b}\\
    \end{aligned}
\end{equation}
The shift operation introduces the possibility of regions outside the original feasible range $[\mathbf{a},\mathbf{b}]$ where, in certain instances, points with lower function values than the global optimum within the original feasible range may exist.
To tackle this issue, the transformed vector $\hat{\mathbf{x}}$ \eqref{eq:shift} is restricted to lie within the range $[\mathbf{a},\mathbf{b}]$ using min-max functions.

Nevertheless, one drawback of this approach is that the functions become ``flat'' in areas where the min-max restriction is applied. 
These flat regions increase in size as the value of $\lambda$ increases.
To address this concern, we opted to limit the range of the randomly generated shift vector by assigning a uniformly distributed random multiplication rate $\rho_j \in [0, 0.1]$ to each dimension $j = 1,...,n$.

For convenient access to all test problems utilized in this paper and to replicate the random shift vectors, we created a dedicated \matlab{} script in the ``Scripts/MPC'' directory of the GitHub repository (\url{https://github.com/blockchain-group/DIRECTGO}). 
These scripts serve as valuable tools for reproducing the findings presented in this investigation and for comparing and evaluating newly developed algorithms.

\subsection{Setup and fundamental basis for algorithm comparison}

All computations were executed on an Intel(R) Core$^\textit{TM}$ i5-10400 @ 2.90GHz Processor running \matlab{} R2023a. The algorithms' solutions were compared with the globally optimal solution for each problem, and we considered the solver successful when the objective function value of a solution was within $0.01 \%$ of the global optimum.

For all analytical test cases with a known global optimum $f^*$, we employed a stopping criterion based on the percent error $(pe)$, as defined below:
\begin{equation}
    \label{eq:pe}
    \ pe = 100 \times
    \begin{cases}
        \frac{f({\mathbf{x}}) - f^*}{\mid f^* \mid},  & f^* \neq 0 \\
        f({\mathbf{x}}),  & f^* = 0
    \end{cases}
\end{equation}
The algorithms were terminated under the following conditions: 
\begin{itemize}
    \item When the $(pe)$ became smaller than $\varepsilon_{\rm pe} = 0.01$.
    \item When the number of function evaluations exceeded the prescribed limit $M_{\rm max} = n \times 10^5$.
    \item When the execution time exceeded $T_{\rm max} = 30$ CPU minutes. In such cases, the final result was set to $n \times 10^5$ to facilitate further processing of the result.
\end{itemize}

\subsection{Algorithm design in \directgen}

Considering that the developed \directgen{} software allows for a large number of combinations, identifying the most effective ones may require a substantial amount of time and effort. 
Therefore, we cannot guarantee that the algorithms presented are the most efficient within \directgen.
Furthermore, the benchmark set includes numerous distinct problems, such as discontinuous, non-differentiable, multi-modal, non-symmetric, and plateau functions. 
It is improbable that a single combination will be the most efficient for all of these diverse problem types.

According to the no-free lunch theorem for optimization~\cite{Wolpert1995NoFL}, there exists no universal optimization algorithm that performs optimally on all types of optimization problems. 
As a result, certain modifications and additions to specific algorithms may not enhance performance on all problems and could even lead to a decline in performance in certain cases.
Therefore, the most optimal approach would involve leveraging machine learning-enhanced automated algorithm selection techniques~\cite{Kerschke20198} to generate algorithms tailored to specific problems. 
However, this avenue remains a part of our future work and has yet to be explored.

To showcase the benefits of \directgen{} software, we conducted an experiment involving five existing \direct-type algorithms: \dtcg~\cite{Stripinis2022wcgo}, \halrect~\cite{Stripinis2022arx}, \mrdirect~\cite{Liu2015}, \birmin~\cite{Paulavicius2019:eswa}, and \dirmin~\cite{Liuzzi2016}. 
Our aim was to improve their average performance across a designated set of test problems by introducing new algorithmic steps or substituting existing ones.

In \Cref{tab:Algorithms}, we present five variants for each of the five selected algorithms, with their improved versions. 
For pure algorithms of \direct-type, which are characterized by slow solution refinement, enhancing their performance was achieved by incorporating local search techniques. 
On the other hand, for hybrid methods, we made different adjustments to improve their performance.
Specifically, for the \birmin{} algorithm, our goal was to increase the number of evaluations per iteration through enhancements, while for the \dirmin{} algorithm, we pursued the opposite approach.

\begin{table}[ht]
    \caption{Description of used parameters for each selected algorithm and their improved versions in \directgen. The \textcolor{blue}{blue} color indicates the parameter that has been substituted or has been added.}
    \label{tab:Algorithms}
    \resizebox{\textwidth}{!}{
    \begin{tabular}{llllll}
        \toprule
        \textbf{Original} algorithm parameters & \dtcg & \halrect & \mrdirect & \birmin & \dirmin \\
        \midrule
        \texttt{Partitioning.Strategy} & $'\texttt{DTC}'$ & $'\texttt{DBC}'$ & $'\texttt{DTC}'$ & $'\texttt{DBDP}'$ & $'\texttt{DTC}'$  \\
        \texttt{Partitioning.SubSides} & $'\texttt{One}'$ & $'\texttt{All}'$ & $'\texttt{All}'$ & $'\texttt{One}'$ & $'\texttt{All}'$ \\
        \texttt{Selection.AggrFuncVal} & $'\texttt{Midpoint}'$ & $'\texttt{MidMin}'$ & $'\texttt{Midpoint}'$ & $'\texttt{Min}'$ & $'\texttt{Midpoint}'$ \\
        \texttt{Selection.CandMeasure} & $'\texttt{Diagonal}'$ & $'\texttt{Diagonal}'$ & $'\texttt{Diagonal}'$ & $'\texttt{Diagonal}'$ & $'\texttt{Diagonal}'$ \\
        \texttt{Selection.Strategy} & $'\texttt{Pareto}'$ & $'\texttt{Aggressive}'$ & $'\texttt{Original}'$ & $'\texttt{Original}'$ & $'\texttt{Original}'$ \\
        \texttt{Selection.EqualCand} & $'\texttt{One}'$ & $'\texttt{One}'$ & $'\texttt{All}'$ & $'\texttt{One}'$ & $'\texttt{All}'$ \\
        \texttt{Selection.SolRefin} & $'\texttt{Off}'$ & $'\texttt{Off}'$ & $'\texttt{Min}'$ & $'\texttt{Min}'$ & $'\texttt{Min}'$ \\
        \texttt{Selection.Ep} & $-$ & $-$ & $0.0001$ & $0.0001$ & $0.0001$ \\
        \texttt{Selection.ControlEp} & $'\texttt{Off}'$ & $'\texttt{Off}'$ & $'\texttt{MultiLevel1}'$ & $'\texttt{Off}'$ & $'\texttt{Off}'$ \\
        \texttt{Selection.GloballyBiased} & $'\texttt{Off}'$ & $'\texttt{Off}'$ & $'\texttt{Off}'$ & $'\texttt{On}'$ & $'\texttt{Off}'$ \\
        \texttt{Selection.TwoPhase} & $'\texttt{On}'$ & $'\texttt{Off}'$ & $'\texttt{Off}'$ & $'\texttt{Off}'$ & $'\texttt{Off}'$ \\
        \texttt{Hybridization.Strategy} & $'\texttt{Off}'$ & $'\texttt{Off}'$ & $'\texttt{Off}'$ & $'\texttt{Single}'$ & $'\texttt{Aggressive}'$ \\
        \texttt{Hybridization.LocalSearch} & $-$ & $-$ & $-$ & $'\texttt{interior-point}'$ & $'\texttt{interior-point}'$ \\
        \texttt{Hybridization.MaxIterations} & $-$ & $-$ & $-$ & $1000$ & $1000$ \\
        \texttt{Hybridization.MaxEvaluations} & $-$ & $-$ & $-$ & $3000$ & $3000$ \\
        \midrule
        \textbf{Improved} algorithm parameters & \dtcg & \halrect & \mrdirect & \birmin & \dirmin \\
        \midrule
        \texttt{Partitioning.Strategy} & $'\texttt{DTC}'$ & $'\texttt{DBC}'$ & $'\texttt{DTC}'$ & $'\texttt{DBDP}'$ & $'\texttt{DTC}'$  \\
        \texttt{Partitioning.SubSides} & $'\texttt{One}'$ & $'\texttt{All}'$ & $'\texttt{All}'$ & $'\textcolor{blue}{\texttt{All}}'$ & $'\texttt{All}'$ \\
        \texttt{Selection.AggrFuncVal} & $'\texttt{Midpoint}'$ & $'\texttt{MidMin}'$ & $'\texttt{Midpoint}'$ & $'\texttt{Min}'$ & $'\texttt{Midpoint}'$ \\
        \texttt{Selection.CandMeasure} & $'\texttt{Diagonal}'$ & $'\texttt{Diagonal}'$ & $'\texttt{Diagonal}'$ & $'\texttt{Diagonal}'$ & $'\textcolor{blue}{\texttt{LongSide}}'$ \\
        \texttt{Selection.Strategy} & $'\texttt{Pareto}'$ & $'\texttt{Aggressive}'$ & $'\texttt{Original}'$ & $'\textcolor{blue}{\texttt{Pareto}}'$ & $'\texttt{Original}'$ \\
        \texttt{Selection.EqualCand} & $'\texttt{One}'$ & $'\texttt{One}'$ & $'\texttt{All}'$ & $'\texttt{One}'$ & $'\textcolor{blue}{\texttt{One}}'$ \\
        \texttt{Selection.SolRefin} & $'\texttt{Off}'$ & $'\texttt{Off}'$ & $'\texttt{Min}'$ & $'\texttt{Min}'$ & $'\textcolor{blue}{\texttt{Median}}'$ \\
        \texttt{Selection.Ep} & $-$ & $-$ & $0.0001$ & $0.0001$ & $0.0001$ \\
        \texttt{Selection.ControlEp} & $'\texttt{Off}'$ & $'\texttt{Off}'$ & $'\texttt{MultiLevel1}'$ & $'\texttt{Off}'$ & $'\texttt{Off}'$ \\
        \texttt{Selection.GloballyBiased} & $'\texttt{Off}'$ & $'\texttt{Off}'$ & $'\texttt{Off}'$ & $'\texttt{On}'$ & $'\texttt{Off}'$ \\
        \texttt{Selection.TwoPhase} & $'\texttt{On}'$ & $'\texttt{Off}'$ & $'\texttt{Off}'$ & $'\texttt{Off}'$ & $'\texttt{Off}'$ \\
        \texttt{Hybridization.Strategy} & $'\textcolor{blue}{\texttt{Single}}'$ & $'\textcolor{blue}{\texttt{Aggressive}}'$ & $'\textcolor{blue}{\texttt{Clustering}}'$ & $'\texttt{Single}'$ & $'\texttt{Aggressive}'$ \\
        \texttt{Hybridization.LocalSearch} & $'\textcolor{blue}{\texttt{sqp}}'$ & $'\textcolor{blue}{\texttt{sqp}}'$ & $'\textcolor{blue}{\texttt{sqp}}'$ & $'\textcolor{blue}{\texttt{sqp}}'$ & $'\texttt{interior-point}'$ \\
        \texttt{Hybridization.MaxIterations} & $\textcolor{blue}{1000}$ & $\textcolor{blue}{1000}$ & $\textcolor{blue}{1000}$ & $1000$ & $1000$ \\
        \texttt{Hybridization.MaxEvaluations} & $\textcolor{blue}{3000}$ & $\textcolor{blue}{3000}$ & $\textcolor{blue}{3000}$ & $3000$ & $3000$ \\
        \bottomrule
    \end{tabular}}
\end{table}

It is essential to note that the construction of the original algorithms within \directgen{} may not always produce identical results to the implementations provided in \texttt{DIRECTGO}~\cite{Stripinis2021c}. 
The discrepancy in the results can be attributed to the numerical tolerances used in the implementations, which play a critical role in the outcome.
For instance, authors might employ rounding on hyper-rectangle measure sizes, enabling them to group extremely small hyper-rectangles together. 
Additionally, they might consider two function values identical if their difference is below a certain threshold. 
These variations in the implementations can significantly impact the selection of POCs.

\subsection{Results and discussions}

In this section, we conduct a performance evaluation of ten \direct-type algorithms, five of which are newly generated with \directgen{}.
The experimental results presented in this evaluation can also be accessed digitally in the ``Results/MPC'' directory of the GitHub repository, available at \url{https://github.com/blockchain-group/DIRECTGO}.

\subsubsection{Comparison of success rates and function evaluations utilization}

\Cref{tab:Results} provides an overview of the success rates achieved by the ten \direct-type approaches considered on various subsets of the \directgolib{} test problems.
In particular, improvements that effectively improve the performance of the original algorithm are highlighted in \textcolor{green}{green}, while those that lead to deteriorating results are marked in \textcolor{red}{red}.
The most remarkable enhancements in success rates were observed in the case of the algorithm that performed worst in this study (\mrdirect) after applying the improvements. 
Its enhanced version yielded a remarkable increase in the success rate of $15.78 \%$.
Moreover, the most significant improvements were evident in the resolution of uni-modal problems, where the pure \mrdirect{} version failed to locate the desired solutions within the allocated evaluation budget in $32.41 \%$ fewer instances.

\begin{table}[ht]
    \caption{Comparison of the success rates of different algorithms in solving test problems with various characteristics.}
    \label{tab:Results}
	\begin{tabular*}{\textwidth}{@{\extracolsep{\fill}}llcccccc}
		\toprule
        \multirow{3}{*}{Algorithm} & \multicolumn{7}{c}{Percentage of solved problems} \\
        \cmidrule{2-8} 
        & \multirow{2}{*}{Overall} & \multicolumn{2}{c}{Separability} & \multicolumn{2}{c}{Multi-modality} & \multicolumn{2}{c}{Scalability} \\
        && $+$ & $-$ & $+$ & $-$ & $+$ & $-$ \\
		\midrule
        Impr. \dtcg & $\textcolor{green}{82.18}$ & $91.71$ & $\textcolor{green}{77.62}$ & $\textcolor{green}{77.30}$ & $\textcolor{green}{98.62}$ & $\textcolor{green}{78.92}$ & $94.12$ \\ 
        Orig. \dtcg & $80.60$ & $91.71$ & $75.29$ & $75.66$ & $97.24$ & $76.91$ & $94.12$ \\
        \midrule
        Impr. \halrect & $\textcolor{green}{78.08}$ & $\textcolor{green}{87.32}$ & $\textcolor{green}{73.66}$ & $\textcolor{green}{72.19}$ & $\textcolor{green}{97.93}$ & $\textcolor{green}{73.69}$ & $\textcolor{green}{94.12}$ \\
        Orig. \halrect & $67.35$ & $76.59$ & $62.94$ & $62.58$ & $83.45$ & $62.45$ & $85.29$ \\
        \midrule
        Impr. \mrdirect & $\textcolor{green}{64.20}$ & $\textcolor{green}{78.05}$ & $\textcolor{green}{57.58}$ & $\textcolor{green}{55.42}$ & $\textcolor{green}{93.79}$ & $\textcolor{green}{59.84}$ & $\textcolor{green}{80.15}$ \\
        Orig. \mrdirect & $48.42$ & $65.85$ & $40.09$ & $44.58$ & $61.38$ & $43.17$ & $67.65$ \\ 
        \midrule
        Impr. \birmin & $\textcolor{green}{75.08}$ & $\textcolor{green}{83.90}$ & $\textcolor{green}{70.86}$ & $\textcolor{green}{68.30}$ & $\textcolor{green}{97.93}$ & $\textcolor{green}{70.68}$ & $91.17$ \\
        Orig. \birmin & $70.66$ & $82.43$ & $65.03$ & $63.60$ & $94.48$ & $65.06$ & $91.17$ \\
        \midrule
        Impr. \dirmin & $\textcolor{red}{76.34}$ & $\textcolor{red}{84.39}$ & $\textcolor{red}{72.49}$ & $\textcolor{red}{70.34}$ & $96.55$ & $\textcolor{red}{71.28}$ & $94.85$ \\
        Orig. \dirmin & $77.76$ & $84.88$ & $74.36$ & $72.19$ & $96.55$ & $73.09$ & $94.85$ \\
	\bottomrule
    \end{tabular*}
\end{table}

Among the pure \direct-type algorithms, the \dtcg{} algorithm exhibited the lowest increase in success rates. 
When considering the allocated budget for function evaluations, the improved algorithm \dtcg{} failed to provide a solution to the $113$ problems, while the original version struggled with the $123$ problems.
An important observation is that the original algorithm \dtcg{} performed quite well, surpassing the overall performance of the improved versions of other less efficient algorithms.

The enhancements in the hybrid algorithms resulted in increased success rates only for the \birmin{} algorithm, whereas the success rates for the \dirmin{} algorithms exhibited a slight deterioration in most of the subsets considered.
Despite the improvement achieved in the \birmin{} algorithm, it still remained outperformed by both versions of the \dirmin{} algorithm in almost all cases.

\begin{figure}[ht]
    \resizebox{\textwidth}{!}{
    \begin{tikzpicture}
        \begin{groupplot}[group style={group size=1 by 5,x descriptions at=edge bottom,y descriptions at=edge left,vertical sep=11pt,horizontal sep=0pt,},]
            \nextgroupplot[
                title={\dtcg},title style={yshift=-2ex}, height=0.215\textwidth,width=\textwidth,
                xmin=1,xmax=100000, xtick distance=10,enlarge x limits=0.01, xmode=log, ytick={1,2},
                yticklabels={Original, Improved},
                yticklabel style = {font=\small},enlarge y limits=0.15,
                xticklabel style = {font=\large},
                grid=both,grid style={line width=.1pt, draw=gray!15},
                boxplot/draw direction=x,boxplot/whisker extend=0.5,
                boxplot/every median/.style={red, thick, line width=1pt},
                boxplot/every whisker/.style={black, ultra thick, dotted, line width=1pt},
                boxplot/every box/.style={thick},boxplot/box extend=0.8,
                every boxplot/.style={mark=+,every mark/.append style={mark size=2pt}, mark options={draw=red}}
                ]
               ]
               \foreach \n in {0,...,1} {
                    \addplot+[boxplot, fill=none, draw=black, solid, line width=.2pt] table[y index=\n] {boxpl_fevals.txt};
                }
            \nextgroupplot[
                title={\halrect},title style={yshift=-2ex}, height=0.215\textwidth,width=\textwidth,
                xmin=1,xmax=100000, xtick distance=10,enlarge x limits=0.01, xmode=log, ytick={1,2},
                yticklabels={Original, Improved},
                yticklabel style = {font=\small},enlarge y limits=0.15,
                xticklabel style = {font=\large},
                grid=both,grid style={line width=.1pt, draw=gray!15},
                boxplot/draw direction=x,boxplot/whisker extend=0.5,
                boxplot/every median/.style={red, thick, line width=1pt},
                boxplot/every whisker/.style={black, ultra thick, dotted, line width=1pt},
                boxplot/every box/.style={thick},boxplot/box extend=0.8,
                every boxplot/.style={mark=+,every mark/.append style={mark size=2pt}, mark options={draw=red}}
                ]
               ]
               \foreach \n in {2,...,3} {
                \addplot+[boxplot, fill=none, draw=black, solid, line width=.2pt] table[y index=\n] {boxpl_fevals.txt};
                }
            \nextgroupplot[
                title={\mrdirect},title style={yshift=-2ex}, height=0.215\textwidth,width=\textwidth,
                xmin=1,xmax=100000, xtick distance=10,enlarge x limits=0.01, xmode=log, ytick={1,2},
                yticklabels={Original, Improved},
                yticklabel style = {font=\small},enlarge y limits=0.15,
                xticklabel style = {font=\large},
                grid=both,grid style={line width=.1pt, draw=gray!15},
                boxplot/draw direction=x,boxplot/whisker extend=0.5,
                boxplot/every median/.style={red, thick, line width=1pt},
                boxplot/every whisker/.style={black, ultra thick, dotted, line width=1pt},
                boxplot/every box/.style={thick},boxplot/box extend=0.8,
                every boxplot/.style={mark=+,every mark/.append style={mark size=2pt}, mark options={draw=red}}
                ]
               ]
               \foreach \n in {4,...,5} {
                \addplot+[boxplot, fill=none, draw=black, solid, line width=.2pt] table[y index=\n] {boxpl_fevals.txt};
                }
            \nextgroupplot[
                title={\birmin},title style={yshift=-2ex}, height=0.215\textwidth,width=\textwidth,
                xmin=1,xmax=100000, xtick distance=10,enlarge x limits=0.01, xmode=log, ytick={1,2},
                yticklabels={Original, Improved},
                yticklabel style = {font=\small},enlarge y limits=0.15,
                xticklabel style = {font=\large},
                grid=both,grid style={line width=.1pt, draw=gray!15},
                boxplot/draw direction=x,boxplot/whisker extend=0.5,
                boxplot/every median/.style={red, thick, line width=1pt},
                boxplot/every whisker/.style={black, ultra thick, dotted, line width=1pt},
                boxplot/every box/.style={thick},boxplot/box extend=0.8,
                every boxplot/.style={mark=+,every mark/.append style={mark size=2pt}, mark options={draw=red}}
                ]
               ]
               \foreach \n in {6,...,7} {
                \addplot+[boxplot, fill=none, draw=black, solid, line width=.2pt] table[y index=\n] {boxpl_fevals.txt};
                }
            \nextgroupplot[
                title={\dirmin},title style={yshift=-2ex}, height=0.215\textwidth,width=\textwidth,
                xmin=1,xmax=100000, xtick distance=10,enlarge x limits=0.01, xmode=log,
                xlabel={Function evaluations / dimension}, ytick={1,2},
                yticklabels={Original, Improved},
                yticklabel style = {font=\small},enlarge y limits=0.15,
                xticklabel style = {font=\large},
                grid=both,grid style={line width=.1pt, draw=gray!15},
                boxplot/draw direction=x,boxplot/whisker extend=0.5,
                boxplot/every median/.style={red, thick, line width=1pt},
                boxplot/every whisker/.style={black, ultra thick, dotted, line width=1pt},
                boxplot/every box/.style={thick},boxplot/box extend=0.8,
                every boxplot/.style={mark=+,every mark/.append style={mark size=2pt}, mark options={draw=red}}
                ]
               ]
               \foreach \n in {8,...,9} {
                \addplot+[boxplot, fill=none, draw=black, solid, line width=.2pt] table[y index=\n] {boxpl_fevals.txt};
                }
        \end{groupplot}
    \end{tikzpicture}}
    \caption{Box plot graphical comparison of algorithms performance based on function evaluations per dimension across all test problems.}
    \label{fig:boxplot1}
\end{figure}
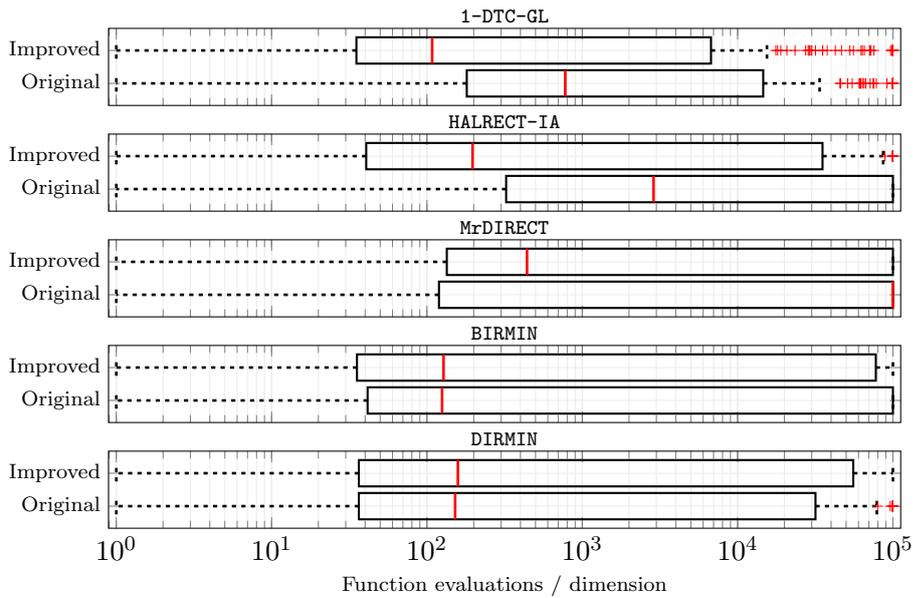

\Cref{fig:boxplot1} presents a box plot that compares algorithms based on function evaluations per dimension on all test problems. 
An important distinction between pure and hybrid algorithms is that pure algorithms generally require more function evaluations, even for relatively simple optimization problems. 
On the contrary, hybrid algorithms demonstrate the ability to solve such problems quickly and efficiently.
Among the algorithms, almost all hybrid algorithms achieved similar lowest first-quartile values, indicating that these methods could solve at least $25\%$ of the test problems faster than pure algorithms.
Specifically, four algorithms (original and improved \dirmin, improved \birmin, and improved \dtcg) were in the lowest first quartile.
On the contrary, the original \halrect{} and original \dtcg{} algorithms exhibited the worst first-quartile performance, each requiring approximately nine and six times more function evaluations, respectively, than the best-performing algorithm, \dirmin.

The improved algorithm \dtcg{} demonstrated the best median value, while its pure counterpart, the original version \dtcg{}, had the third worst median value in these studies. 
The addition of the local search procedure to the \dtcg{} algorithm reduced the median value by nearly eight times, resulting in a significant improvement in its performance.
Interestingly, the median value of the original algorithm \mrdirect{} is equal to the maximum number of function evaluation budgets ($M_{\rm max}$), indicating that the algorithm could not solve more than half of the test problems. 
However, its improved version exhibited a significantly higher median value.
When comparing the third-quartile values, four algorithms reached the $M_{\rm max}$ value in the third quartile, suggesting that these algorithms could not solve more than $25\%$ of the problems. 
Only six algorithms achieved values lower than the maximum evaluation budget. 
Among these, the improved algorithm \dtcg{} achieved the lowest third-quartile value, approximately half that of the second-best pure algorithm, the original algorithm \dtcg.

\subsubsection{Analysis of results across different subsets of problems}

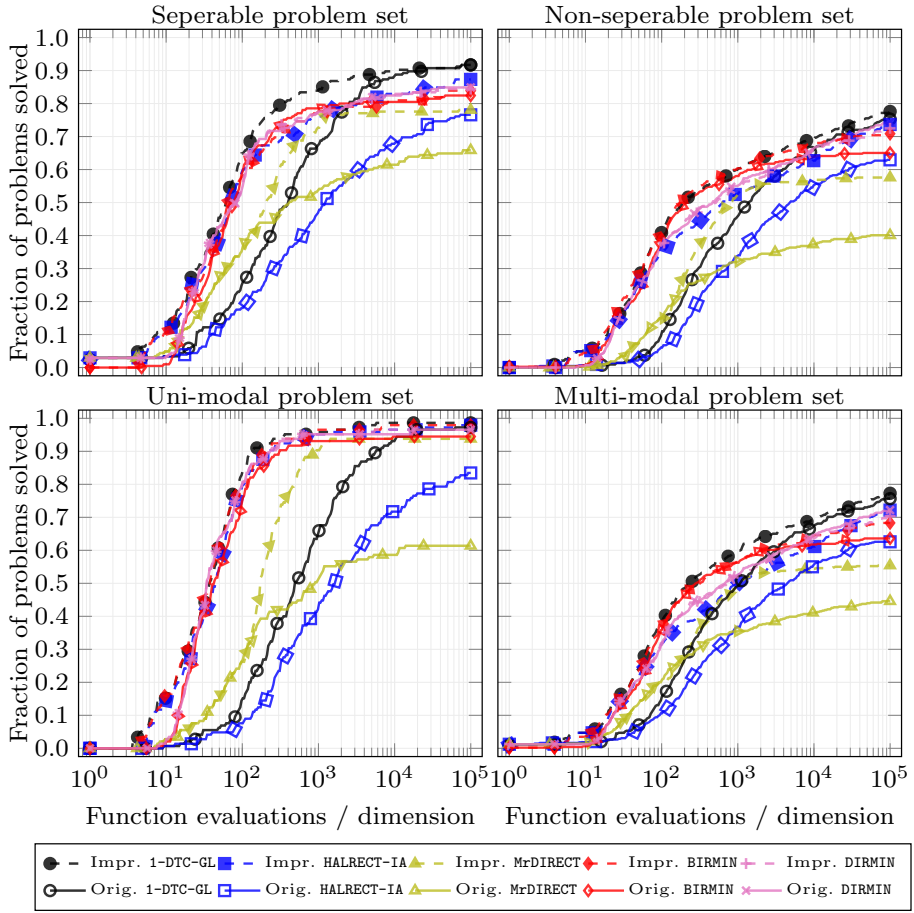
\begin{figure}
    \resizebox{\textwidth}{!}{
    \begin{tikzpicture}
        \begin{groupplot}[group style={group size=2 by 2,x descriptions at=edge bottom,y descriptions at=edge left,vertical sep=11pt,horizontal sep=5pt,},]
            \nextgroupplot[
            title={Seperable problem set},title style={yshift=-2.5ex},
            ymin=0,ymax=1,ytick distance=0.1,enlarge y limits=0.025,
            ylabel={Fraction of problems solved},ylabel shift = -5pt,
            yticklabels={0.0,0.1,0.2,0.3,0.4,0.5,0.6,0.7,0.8,0.9,1.0},ytick={0.0,0.1,0.2,0.3,0.4,0.5,0.6,0.7,0.8,0.9,1.0},
            yticklabel style = {font=\small},
            xmin=1,xmax=100000,xtick distance=10,enlarge x limits=0.03,xmode=log,
            grid=both,grid style={line width=.1pt, draw=gray!15},
            tick label style={font=\small},height=0.45\textwidth,width=0.5\textwidth,
            ]
            \addplot[s_dc]	table[x=T,y=b1]	{overall.txt};
            \addplot[s_dci]	table[x=T,y=b2]	{overall.txt};
            \addplot[s_hl]	table[x=T,y=b3]	{overall.txt};
            \addplot[s_hli]	table[x=T,y=b4]	{overall.txt};
            \addplot[s_mr]	table[x=T,y=b5] {overall.txt};
            \addplot[s_mri]	table[x=T,y=b6]	{overall.txt};
            \addplot[s_bm]	table[x=T,y=b7]	{overall.txt};
            \addplot[s_bmi]	table[x=T,y=b8]	{overall.txt};
            \addplot[s_dm]	table[x=T,y=b9]	{overall.txt};
            \addplot[s_dmi]	table[x=T,y=b10]{overall.txt};
            \nextgroupplot[
            title={Non-seperable problem set},title style={yshift=-2.5ex},
            ymin=0,ymax=1,ytick distance=0.1,enlarge y limits=0.025,
            xmin=1,xmax=100000,xtick distance=10,enlarge x limits=0.03,xmode=log,
            grid=both,grid style={line width=.1pt, draw=gray!15},
            tick label style={font=\normalsize},height=0.45\textwidth,width=0.5\textwidth,
            ]
            \addplot[s_dc]	table[x=T,y=c1]	{overall.txt};
            \addplot[s_hl]	table[x=T,y=c3]	{overall.txt};
            \addplot[s_mr]	table[x=T,y=c5] {overall.txt};
            \addplot[s_bm]	table[x=T,y=c7]	{overall.txt};
            \addplot[s_dm]	table[x=T,y=c9]	{overall.txt};
            \addplot[s_dci]	table[x=T,y=c2]	{overall.txt};
            \addplot[s_hli]	table[x=T,y=c4]	{overall.txt};
            \addplot[s_mri]	table[x=T,y=c6]	{overall.txt};
            \addplot[s_bmi]	table[x=T,y=c8]	{overall.txt};
            \addplot[s_dmi]	table[x=T,y=c10]{overall.txt};
            \nextgroupplot[
            title={Uni-modal problem set},title style={yshift=-2.5ex},
            ymin=0,ymax=1,ytick distance=0.1,enlarge y limits=0.025,
            ylabel={Fraction of problems solved},ylabel shift = -5pt,
            xlabel={Function evaluations / dimension},
            yticklabels={0.0,0.1,0.2,0.3,0.4,0.5,0.6,0.7,0.8,0.9,1.0},ytick={0.0,0.1,0.2,0.3,0.4,0.5,0.6,0.7,0.8,0.9,1.0},
            xticklabel style = {font=\small},
            yticklabel style = {font=\small},
            xmin=1,xmax=100000,xtick distance=10,enlarge x limits=0.03,xmode=log,
            grid=both,grid style={line width=.1pt, draw=gray!15},
            tick label style={font=\small},height=0.45\textwidth,width=0.5\textwidth,
            ]
            \addplot[s_dci]	table[x=T,y=d2]	{overall.txt};
            \addplot[s_hli]	table[x=T,y=d4]	{overall.txt};
            \addplot[s_mri]	table[x=T,y=d6]	{overall.txt};
            \addplot[s_bmi]	table[x=T,y=d8]	{overall.txt};
            \addplot[s_dmi]	table[x=T,y=d10]{overall.txt};
            \addplot[s_dc]	table[x=T,y=d1]	{overall.txt};
            \addplot[s_hl]	table[x=T,y=d3]	{overall.txt};
            \addplot[s_mr]	table[x=T,y=d5] {overall.txt};
            \addplot[s_bm]	table[x=T,y=d7]	{overall.txt};
            \addplot[s_dm]	table[x=T,y=d9]	{overall.txt};
            \nextgroupplot[
            title={Multi-modal problem set},title style={yshift=-2.5ex},
            ymin=0,ymax=1,ytick distance=0.1,enlarge y limits=0.025,
            xlabel={Function evaluations / dimension},
            xmin=1,xmax=100000,xtick distance=10,enlarge x limits=0.03,xmode=log,
            xticklabel style = {font=\small},
            grid=both,grid style={line width=.1pt, draw=gray!15},
            tick label style={font=\small},height=0.45\textwidth,width=0.5\textwidth,
            legend columns=5,legend style={at={(1.025,-0.25)},font=\tiny, column sep=-1pt,legend cell align={left},}
            ]
            \addplot[s_dci]	table[x=T,y=e2]	{overall.txt};
            \addplot[s_hli]	table[x=T,y=e4]	{overall.txt};
            \addplot[s_mri]	table[x=T,y=e6]	{overall.txt};
            \addplot[s_bmi]	table[x=T,y=e8]	{overall.txt};
            \addplot[s_dmi]	table[x=T,y=e10]{overall.txt};       
            \addplot[s_dc]	table[x=T,y=e1]	{overall.txt};
            \addplot[s_hl]	table[x=T,y=e3]	{overall.txt};
            \addplot[s_mr]	table[x=T,y=e5] {overall.txt};
            \addplot[s_bm]	table[x=T,y=e7]	{overall.txt};
            \addplot[s_dm]	table[x=T,y=e9]	{overall.txt};
            \legend{Impr. \dtcg, Impr. \halrect, Impr. \mrdirect, Impr. \birmin, Impr. \dirmin, Orig. \dtcg, Orig. \halrect, Orig. \mrdirect, Orig. \birmin, Orig. \dirmin} 
        \end{groupplot}
    \end{tikzpicture}}
    \caption{Data profiles: The horizontal axis represents the number of function evaluations per dimension, while the vertical axis represents the fraction of solved problems}
    \label{figs_all}
\end{figure}

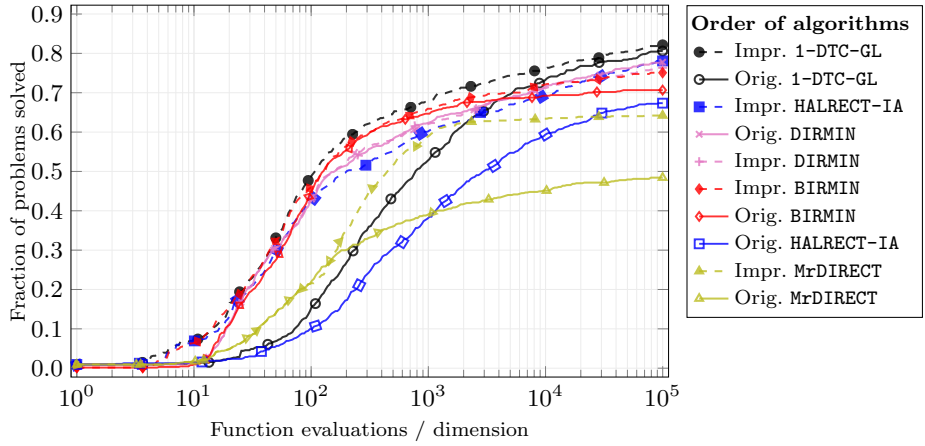
\begin{figure}[ht]
\resizebox{\textwidth}{!}{
    \begin{tikzpicture}
        \begin{axis}[
            ymin=0,ymax=0.9,ytick distance=0.1,enlarge y limits=0.025,
            ylabel={Fraction of problems solved},ylabel shift = -5pt,
            xlabel={Function evaluations / dimension},
            yticklabels={0.0,0.1,0.2,0.3,0.4,0.5,0.6,0.7,0.8,0.9,1.0},ytick={0.0,0.1,0.2,0.3,0.4,0.5,0.6,0.7,0.8,0.9,1.0},
            xmin=1,xmax=100000,xtick distance=10,enlarge x limits=0.01,xmode=log,
            grid=both,grid style={line width=.1pt, draw=gray!15},
            tick label style={font=\normalsize},height=0.55\textwidth,width=0.8\textwidth,
            legend pos=outer north east,legend style={row sep=0pt,font=\small},legend cell align={left},
            ]
            \addlegendimage{empty legend}
            \addplot[s_dci]	table[x=T,y=a2]	{overall.txt};
            \addplot[s_dc]	table[x=T,y=a1]	{overall.txt};
            \addplot[s_hli]	table[x=T,y=a4]	{overall.txt};
            \addplot[s_dm]	table[x=T,y=a9]	{overall.txt};
            \addplot[s_dmi]	table[x=T,y=a10] {overall.txt};
            \addplot[s_bmi]	table[x=T,y=a8]	{overall.txt};
            \addplot[s_bm]	table[x=T,y=a7]	{overall.txt};
            \addplot[s_hl]	table[x=T,y=a3]	{overall.txt};
            \addplot[s_mri]	table[x=T,y=a6]	{overall.txt};
            \addplot[s_mr]	table[x=T,y=a5]	{overall.txt};
            \legend{\hspace{-.6cm}\textbf{Order of algorithms}, Impr. \dtcg, Orig. \dtcg, Impr. \halrect, Orig. \dirmin, Impr. \dirmin, Impr. \birmin, Orig. \birmin, Orig. \halrect, Impr. \mrdirect, Orig. \mrdirect} 
        \end{axis}
    \end{tikzpicture}}
    \caption{Data profiles: The horizontal axis represents the number of function evaluations per dimension, while the vertical axis represents the fraction of solved problems \label{fig1}}
\end{figure}

The data profiles~\cite{More2009} depicted in \Cref{figs_all} showcase how all algorithms perform on test problems with various properties of \directgolib. 
These profiles provide a comprehensive view of algorithm performance across different types of problems.
Meanwhile, the data profiles in \Cref{fig1} offer an overall ranking of the algorithms on all test problems, providing a more focused perspective on their performance in a broader context.

Hybridization of pure \direct-type algorithms significantly impacts the results, particularly when dealing with straightforward uni-modal or separable test problems. 
The inclusion of a local search procedure proves to be particularly advantageous for uni-modal problems, as it accelerates the convergence speed to reach optimal solutions more efficiently. 
On the other hand, pure \direct-type algorithms might prioritize the global search and exhaust the evaluation budget without locating the solution within the prescribed accuracy.
As a result, the curves of the improved versions of \dtcg, \halrect, and \mrdirect{} demonstrate significantly better performance than the original versions, especially for small evaluation budgets (${ \leq} 1000{\times}n$).
However, it is worth noting that the most successful pure \direct-type algorithm, \dtcg, eventually achieves nearly identical performance within the maximum evaluation budget, regardless of whether the problems are separable or uni-modal.

The improved hybrid algorithm \birmin{} exhibits slightly lower performance within a small evaluation budget ($M_{\rm max}{\leq} n{\times}10^2$). 
However, as the evaluation budget increases ($M_{\rm max}{\geq} n{\times}10^4$), the improved version outperforms the original version. 
This difference in performance becomes particularly evident when the algorithm is applied to non-separable or multi-modal test problems.

On the other hand, the curves of the two versions of the hybrid algorithm \dirmin{} are almost indistinguishable within a smaller evaluation budget ($M_{\rm max}{\leq} 2n{\times}10^4$). 
However, within a larger evaluation budget, the original algorithm \dirmin{} exhibits slightly better performance.

Based on the four graphs in \Cref{figs_all} and the overall ranking of the algorithms in \Cref{fig1}, a consistent conclusion can be drawn: the performance of the improved and original \dtcg{} algorithm is the most efficient or at least comparable to the best-performing algorithm.
Analyzing the curves, it is evident that the improved algorithm \dtcg{} exhibits the highest efficiency rates across all graphs compared to other algorithms within any evaluation budget in $[0, n{\times}10^5]$. 
Although the original performance of \dtcg{} becomes competitive, it requires a significant budget for the evaluations of functions ($M_{\rm max} {\geq} n{\times}10^4$).
Overall, the improved performance of the \dtcg{} algorithm remains competitive, requiring fewer function evaluations to achieve the desired optimal value in most test functions.

\subsubsection{Statistical analysis of the results}

To validate the results and comparisons between algorithms, as well as to evaluate the significance of improvements achieved by \directgen{}, we conducted the Friedman mean rank test~\cite{friedman1937use} and the non-parametric Wilcoxon signed test~\cite{hollander1999nonparametric} at a significance level of $5 \%$. 
A $p$-value greater than $0.05$ indicates that the difference in results between methods is statistically insignificant.

\Cref{tab:Resultsfriedman} displays the Friedman mean rank values for considered algorithms, utilizing the founded solution values within different evaluation budgets for all test problems. 
The results reveal that the improved versions consistently outperform their original counterparts in all budgets, except for the algorithm \dirmin{}, where the original version obtained a higher rank in one specific evaluation budget ($M_{\rm max} {=} n{\times}10^5$). 
The improvements made to the algorithms have resulted in performance gains ranging from small to significant, as indicated by the higher mean rank values.

\begin{table}[ht]
    \caption{Friedmann mean rank values with different objective function evaluation budgets.}
    \label{tab:Resultsfriedman}
    \begin{tabular*}{\textwidth}{@{\extracolsep{\fill}}lllll}
        \toprule
        \multirow{2}{*}{Algorithm} & \multicolumn{4}{c}{Function evaluation budget $(M_{\rm max})$} \\
        \cmidrule{2-5} 
        & $n{\times}10^2$ & $n{\times}10^3$ & $n{\times}10^4$ & $n{\times}10^5$ \\
        \midrule
        Impr. \dtcg & $4.7610$ & $4.6128$ & $4.6128$ & $4.6601$ \\
        Orig. \dtcg & $6.0095$ & $5.3375$ & $5.3375$ & $4.7831$ \\
        \midrule
        Impr. \halrect & $5.5670$ & $5.5229$ & $5.5229$ & $5.3249$ \\
        Orig. \halrect & $7.2752$ & $7.1372$ & $7.1372$ & $6.1672$ \\
        \midrule
        Impr. \mrdirect & $5.3730$ & $5.2500$ & $5.2500$ & $5.8028$ \\
        Orig. \mrdirect & $5.3730$ & $6.8099$ & $6.8099$ & $7.1435$ \\
        \midrule
        Impr. \birmin & $4.6151$ & $4.6562$ & $4.6562$ & $5.1333$ \\
        Orig. \birmin & $5.0103$ & $5.0765$ & $5.0765$ & $5.6356$ \\
        \midrule
        Impr. \dirmin & $5.4219$ & $5.2886$ & $5.2886$ & $5.2492$ \\
        Orig. \dirmin & $5.5938$ & $5.3084$ & $5.3084$ & $5.1002$ \\
    \bottomrule
    \end{tabular*}
\end{table}

\Cref{tab:Resultswilc} presents the $p$-values obtained by comparing the improved algorithms with their original counterparts, using the solutions found within four evaluation budgets on all test problems.
For the improved algorithm \dtcg{}, there is strong statistical evidence that the improved version significantly outperforms the original version within a small evaluation budget ($M_{\rm max}{\leq} n{\times}10^3$). 
However, as the evaluation budget increases ($M_{\rm max}{>} n{\times}10^3$), the higher $p$-values suggest that the significance of the improvement decreases and the difference between the improved and original versions becomes less statistically significant.
In contrast, the situation is different for the other two pure \direct-type algorithms. 
The improved versions of \halrect{} and \mrdirect{} show no significant improvement compared to the original versions at $M_{\rm max} {=} n{\times}10^2$. 
However, for larger evaluation budgets, the $p$-values are low, indicating that the improvements are statistically significant.

For the \birmin{} algorithm, the $p$-values are low, indicating that the improvement of the improved version of \birmin{} compared to the original version is statistically significant in all these budgets.
Regarding \dirmin{] algorithm, we can conclude that the improved version of the algorithm is statistically better if the evaluation budgets are $M_{\rm max} {=} n{\times}10^2$ and $M_{\rm max} {=} n{\times}10^4$.

\begin{table}[ht]
    \caption{Wilcoxon signed test $p$-values at $5 \%$ significance level, comparing improved vs. original algorithms across various objective function evaluation budgets.}
    \label{tab:Resultswilc}
	\begin{tabular*}{\textwidth}{@{\extracolsep{\fill}}lllll}
		\toprule
        \multirow{2}{*}{Algorithm} & \multicolumn{4}{c}{Function evaluation budget $(M_{\rm max})$} \\
        \cmidrule{2-5} 
        & $n{\times}10^2$ & $n{\times}10^3$ & $n{\times}10^4$ & $n{\times}10^5$ \\
		\midrule
        \dtcg & $6.3325{\times}10^{-3}$ & $2.4241{\times}10^{-8}$ & $2.6994{\times}10^{-1}$ & $5.8529{\times}10^{-1}$ \\
        \halrect & $2.3254{\times}10^{-1}$ & $4.0565{\times}10^{-13}$ & $5.3640{\times}10^{-10}$ & $3.3265{\times}10^{-13}$ \\
        \mrdirect & $1.0000{\times}10^{0}$ & $2.0141{\times}10^{-34}$ & $4.5374{\times}10^{-40}$ & $5.9712{\times}10^{-34}$ \\
        \birmin & $2.5863{\times}10^{-2}$ & $2.8668{\times}10^{-9}$ & $1.0548{\times}10^{-12}$ & $1.5552{\times}10^{-15}$ \\
        \dirmin & $3.0296{\times}10^{-8}$ & $4.6683{\times}10^{-1}$ & $8.2874{\times}10^{-3}$ & $4.4600{\times}10^{-4}$ \\
	\bottomrule
    \end{tabular*}
\end{table}

\section{Conclusions and future works}
\label{sec:conclusion}

This study introduces a novel generalized \direct-type algorithmic framework known as \directgen, for derivative-free global optimization.
The proposed framework empowers users to construct a wide range of \direct-type algorithms. 
Such innovative work can foster the development of new \direct-type algorithms and help identify the most suitable algorithm for various practical applications.

To demonstrate the efficiency of \directgen{}, we enhanced five selected \direct-type algorithms with the goal of improving their performance and solving global optimization problems more effectively. 
Evaluation of these constructed algorithms was carried out using benchmark test functions from \directgolib{}. 
The results were analyzed both graphically and statistically to gain insight into the algorithms' performance. 
The findings concluded that the newly developed versions of the \direct-type algorithms significantly outperformed their original counterparts in most cases.

In conclusion, this paper has focused on box-constrained global optimization problems, but the generalized \direct-type algorithmic framework (\directgen) could potentially be extended to handle constrained cases as well. 
Furthermore, due to the numerous combinations of algorithms within \directgen, manually testing all of them becomes impractical.
Therefore, future research should explore the automation of these processes using advanced machine-learning techniques. 
By automating the algorithmic components process, optimization can become more efficient and effective.


\section*{Data statement}
\textbf{DIRECTGOLib} - \textbf{DIRECT} \textbf{G}lobal \textbf{O}ptimization test problems \textbf{Lib}rary is designed as a continuously-growing open-source GitHub repository to which anyone can easily contribute.
The exact data underlying this article from \directgolib{} can be accessed on GitHub:
\begin{itemize}
    \item \url{https://github.com/blockchain-group/DIRECTGOLib},
\end{itemize}
and used under the MIT license.
We welcome contributions and corrections to it.

\bibliographystyle{spmpsci}      
\bibliography{library}   

\end{document}